\documentclass[11pt,reqno]{amsart}

\usepackage[utf8]{inputenc}
\usepackage[english]{babel}
\usepackage{comment}

\usepackage{mathtools}
\usepackage[shortlabels]{enumitem}

\usepackage{amsthm}
\usepackage{amssymb}
\usepackage{amsfonts}
\usepackage{dsfont}
\usepackage{leftindex}
\usepackage[mathcal]{eucal}
\usepackage{mathrsfs}
\usepackage{svrsymbols}
\usepackage{stix}
\usepackage{xcolor}
\usepackage{times}

\theoremstyle{definition}
\newtheorem{theorem}{Theorem}
\newtheorem{proposition}[theorem]{Proposition}
\newtheorem{lemma}[theorem]{Lemma}
\newtheorem{corollary}[theorem]{Corollary}
\newtheorem{definition}[theorem]{Definition}

\newtheorem{example}[theorem]{Example}

\newtheorem{remark}[theorem]{Remark}
\newtheorem*{remark*}{Remark}

\newtheorem*{claim*}{Claim}
\newtheorem{question}[theorem]{Question}
\newtheorem{maintheorem}{Theorem}

\newcommand{\MM}{\mathbb{M}}
\newcommand{\NN}{\mathbb{N}}
\newcommand{\QQ}{\mathbb{Q}}
\newcommand{\RR}{\mathbb{R}}

\newcommand{\ZZ}{\mathbb{Z}}

\newcommand{\frag}{\mathcal{F}}

\newcommand{\cA}{\mathcal{A}}

\newcommand{\cC}{\mathcal{C}}

\newcommand{\cF}{\mathcal{F}}
\newcommand{\cG}{\mathcal{G}}

\newcommand{\cL}{\mathcal{L}}

\newcommand{\fm}{\mathfrak{m}}
\newcommand{\fa}{\mathfrak{a}}

\newcommand{\fM}{\mathfrak{M}}
\newcommand{\fN}{\mathfrak{N}}
\newcommand{\fP}{\mathfrak{P}}

\DeclareMathOperator{\tp}{tp}
\DeclareMathOperator{\free}{free}
\DeclareMathOperator{\sub}{sub}

\newcommand{\Phitp}{\Phi\text{-}\!\tp}
\newcommand{\fulltp}{\Soo^V\text{-}\!\tp}

\newcommand{\open}[1]{\left\langle\!\left\langle {#1} \right\rangle\!\right\rangle}
\newcommand{\Loo}{L_{\omega_1\omega}}
\newcommand{\Soo}{S_{\omega_1\omega}}

\newcommand{\Surj}{\operatorname{Surj}}
\newcommand{\Xs}{X_V}

\newcommand{\HFC}{\mathcal{H}^{\cF}_C}
\newcommand{\SF}{S^{\cF}}
\newcommand{\SFC}{S^{\cF}_C}

\definecolor{mypink}{RGB}{190, 80, 150}
\newcommand{\synt}[1]{\textcolor{mypink}{#1}}

\definecolor{mypurp}{RGB}{120, 0, 120}
\newcommand{\symb}[1]{\textcolor{mypurp}{#1}}

\newcommand{\Sneg}{\mathop{\synt{s\cdot}}}
\newcommand{\Sor}{\mathop{\synt{\vee}}}
\newcommand{\Splus}{\mathbin{\synt{+}}}
\newcommand{\Sand}{\mathop{\synt{\wedge}}}

\newcommand{\supx}{\mathop{\synt{\sup_x}}}
\newcommand{\supy}{\mathop{\synt{\sup_y}}}
\newcommand{\join}{\mathop{\synt{\bigvee}}}
\newcommand{\meet}{\mathop{\synt{\bigwedge}}}
\newcommand{\Sinf}[1]{\mathop{\synt{\inf_{#1}}}}
\newcommand{\Sdot}{\mathbin{\synt{\dot{-}}}}

\newcommand{\Sx}{\synt{x}}
\newcommand{\Sy}{\synt{y}}

\newcommand{\SM}{\symb{M}}
\newcommand{\SN}{\symb{N}}
\newcommand{\Sf}{\symb{f}}
\newcommand{\Sp}{\symb{p}}
\newcommand{\Sc}{\symb{c}}
\newcommand{\unity}{\symb{1}}
\newcommand{\sortx}{{}_{\Sx}\SM}
\newcommand{\sorty}{{}_{\Sy}\SM}

\renewcommand{\phi}{\varphi}

\title[Model existence and Scott's theorem for infinitary real-valued logic]{Model Existence and Scott's Theorem\\
		for Infinitary Real-Valued Logic:\\
		A Topological Perspective}
\author{Olivia Aubone}
\address{Department of Mathematics, The University of Texas at San Antonio, San Antonio, TX 78249, USA}
\email{olivia.aubone@utsa.edu}

\author{Eduardo Due\~nez}
\address{Department of Mathematics, The University of Texas at San Antonio, San Antonio, TX 78249, USA}
\email{eduardo.duenez@utsa.edu}

\author{Jos\'e Iovino}
\address{Department of Mathematics, The University of Texas at San Antonio, San Antonio, TX 78249, USA}
\email{jose.iovino@utsa.edu}
\thanks{This version: \today. The most recent version of this paper is available at \url{https://utsa-quantum-math.github.io/resources/omitting-types-real-valued-logic/}; it may supersede the version deposited on arXiv.}
\date{}
\dedicatory{To Xavier Caicedo}


\usepackage{hyperref}
\begin{document}

	\begin{abstract}
For countable vocabularies, we extend the omitting types theorem and Scott's theorem to the infinitary logic of Keisler's general real-valued structures, without equality, a prescribed metric or moduli of continuity.
Each formula has its own finite bound; no common bound for the whole language is required.
For every countable fragment, the Henkin types form a Polish, hence \v{C}ech-complete, space, and under the hypothesis of the omitting types theorem, the Henkin types that satisfy the theory and omit a prescribed countable family of types form a dense $G_\delta$ subset of the space of Henkin types that satisfy the theory.
A single Scott sentence characterizes every structure with countable sorts, among all structures, up to approximate back-and-forth equivalence, or, equivalently, up to agreement on the values of all sentences.
On the space of structures with universe the natural numbers, surjections acting on atomic diagrams replace permutations: every approximate-equivalence class is Borel and is separated from the other classes by a Borel Lipschitz function.
The bounded Borel functions invariant under exact equivalence are precisely the values of sentences of the logic enlarged by continuous connectives.
For finite predicate ranges with a uniform normalized gap, the representation holds in the original syntax, yielding López-Escobar's theorem for classical infinitary logic without equality.
\end{abstract}

\maketitle
\enlargethispage{4pt}

\section{Introduction}
\label{sec:intro}

In this paper we extend two theorems of classical infinitary logic, the omitting types theorem and Scott's theorem, to the infinitary extension of Keisler's model theory for real-valued structures~\cite{Keisler:2023}, and we give each a topological form.
In Keisler's general structures, function symbols are interpreted arbitrarily and predicate symbols by bounded real-valued functions with no moduli of continuity; the formulas need not have a common bound, but each carries a bound of its own, and neither logical equality nor a prescribed metric is required.
In Part~I, the omitting types theorem becomes a Baire category argument in the \v{C}ech-complete space of Henkin types of a countable fragment.
In Part~II, Scott's theorem becomes a statement about the Polish space of structures on the natural numbers, where surjections, acting on atomic diagrams, take the place of permutations.

In Part~I we prove model existence and the omitting types theorem for countable fragments of the infinitary logic of general structures, in the generality just described (Theorems~\ref{thm:A} and~\ref{thm:B}).

For countable fragments of classical $\Loo$, the omitting types theorem is proved in Keisler's book by the model existence theorem~\cite[Chapters~3 and~11]{Keisler:1971}, after Makkai had extended Smullyan's method of consistency properties to the infinitary logics on admissible sets~\cite{Makkai:1969}.
Keisler later proved a general omitting types theorem by forcing with finite sets of sentences~\cite{Keisler:1973}, and Morley derived the theorem for countable fragments from the Baire category theorem~\cite{Morley:1974}, so that the omission of countably many types became the density of a countable intersection of open dense sets.
All of these theorems presuppose equality.
Keisler and Miller carried the omitting types theorem into first-order logic without equality, building the model from closed terms and testing the types to be omitted through tuples of terms~\cite[Section~4]{Keisler-Miller:2001}; their setting is two-valued and finitary, and we follow them in building models from closed terms.

In first-order logic, the axioms of equality assert its reflexivity and the interchangeability of equal elements in every atomic formula, and thereby make equality a congruence for the function and predicate symbols.
In real-valued logic, their counterparts are axioms of continuity.\footnote{In \L{}ukasiewicz predicate logic, the congruence axioms for a similarity, or graded equality, make every predicate $1$-Lipschitz for the distance given by one minus the similarity (see \cite[Section~5.6]{Hajek:1998} and the earlier work collected there).}
The role of equality passed to a uniform structure with moduli prescribed by the signature, in Iovino's work on stable Banach space structures~\cite{Iovino:1994,Iovino:1999a,Iovino:1999b}, in Henson and Iovino's model theory of normed space structures~\cite[Sections~2 and~3]{Henson-Iovino:2002}, and in the continuous logic of Ben Yaacov and Usvyatsov~\cite[Section~2]{BenYaacov-Usvyatsov:2010}.
For a countable signature, the uniform structure has a countable base, and it is therefore induced by a single pseudometric (see \cite[Section~8.1]{Engelking:1989}).
In that setting Eagle, drawing in part on work of Caicedo and Iovino~\cite{Caicedo-Iovino:2014}, applied \v{C}ech's form of the Baire category theorem to a space of countable metric structures~\cite{Eagle:2014}.
His theorem is formulated for structures whose signature prescribes a master metric, the real-valued counterpart of equality, together with a modulus of uniform continuity for every symbol, and whose predicates take values in $[0,1]$~\cite[Definitions~1.4 and~1.5]{Eagle:2014}; formulas that are not uniformly bounded, predicates without moduli of continuity, and structures without a prescribed master metric lie outside that setting.
Keisler then obtained the pseudometric from the formulas themselves~\cite[Section~3.3]{Keisler:2023}, and showed that general structures stand to first-order structures without equality as metric structures with a distance predicate stand to first-order structures with equality~\cite[Section~1]{Keisler:2023}; his model theory is first-order.
Theorems~\ref{thm:A} and~\ref{thm:B} hold in a logic that none of these settings contains, infinitary, real-valued and free of each of the assumptions on which these settings rest.

The infinitary language $\Loo$ extends Keisler's first-order language by countable joins and infima, each sentence taking its values in a bounded interval of its own, so that the formal types, families of such values, are the points of a product of compact intervals.
We expand the vocabulary by countably many new constants of each sort and build models from closed terms, which themselves serve as elements; since the vocabulary has no equality symbol, no quotient is taken.
The Henkin types, defined by the real-valued form of Henkin's requirement of witnesses, are exactly the types of such term structures.
The key observation of Part~I is that, for a countable vocabulary and a countable fragment, the space of Henkin types of the fragment is a $G_\delta$ subset of a compact space, hence Polish and \v{C}ech-complete, so that Morley's method operates in that space.

A type of the fragment in a finite tuple $\bar{x}$ of variables is a set of expressions $\synt{\phi(\bar{x})\geq r}$, with $r$ rational and with $\phi$ a formula of the fragment whose free variables are among $\bar{x}$.
The type is omitted by a Henkin type of the fragment that satisfies none of the instances obtained by substituting closed terms of the expanded vocabulary for these variables.

\begin{maintheorem}
\label{thm:A}
Let a countable vocabulary be given, together with a formal type for a set of sentences of $\Loo$ for this vocabulary.
The formal type is realized by a structure whose sorts are countable if and only if it extends to a Henkin type.
If the given sentences belong to a countable fragment, then the formal type is realized by some structure if and only if it extends to a Henkin type of the fragment.
\end{maintheorem}

\begin{maintheorem}
\label{thm:B}
Let a countable fragment of $\Loo$ for a countable vocabulary be given, together with a theory on sentences of the fragment that has a model and a sequence of types of the fragment.
Suppose that, for every instance of every type of the sequence, the Henkin types of the fragment that satisfy the theory and the instance form a set with empty interior in the space of Henkin types of the fragment that satisfy the theory.
Then some Henkin type of the fragment that satisfies the theory omits every type of the sequence, and the theory has a model whose sorts are countable and in which no tuple of elements satisfies any type of the sequence.
\end{maintheorem}

Theorem~\ref{thm:A} is proved in Section~\ref{sec:Henkin-type} (Proposition~\ref{prop:henkin-sat}), and Theorem~\ref{thm:B} in Section~\ref{sec:omit-type} (Theorem~\ref{thm:omitting-types}), after the \v{C}ech completeness of the space of Henkin types is established in Section~\ref{sec:Henkin-topology} (Corollary~\ref{thm:cech-complete}).
The conclusion is generic in the Henkin space: under the hypothesis of Theorem~\ref{thm:B}, the Henkin types satisfying the theory and omitting every type of the sequence form a dense $G_\delta$ subset of the space of Henkin types satisfying the theory.
Omission can therefore be imposed together with any finite collection of strict rational conditions on sentences of the fragment whose corresponding basic open set in that space is nonempty.
The countable model is obtained from a point of this dense set by the term construction.

The hypothesis of Theorem~\ref{thm:B} is the topological form of the hypothesis of Keisler's general omitting types theorem~\cite[p.~108, General Omitting Types Theorem]{Keisler:1973}.
In classical two-valued logic, with Keisler's $\forall\exists$ formulas taken over the set of all formulas of the fragment, a model omits a type $\Gamma(\bar{x})$ exactly when it satisfies the sentence $\forall\bar{x}\bigvee_{\gamma\in\Gamma}\neg\gamma(\bar{x})$.
For this sentence Keisler asks that every finite set of sentences satisfiable in the class remain satisfiable after the addition of $\bigvee_{\gamma\in\Gamma}\neg\gamma(\bar{c})$, for every tuple $\bar{c}$ of new constants.
In topological terms, every instance has empty interior with respect to the basic sets given by finite sets of sentences.
Since every $\forall\exists$ sentence $\forall\bar{x}\bigvee_{i}\theta_i(\bar{x})$ over the fragment holds in exactly the models that omit the type $\{\neg\theta_i\}_i$, this reading covers every sentence to which Keisler's theorem applies.

Theorem~\ref{thm:B} carries this case of the general omitting types theorem to Keisler's general structures with three refinements: its instances are taken at tuples of closed terms, as in the work of Keisler and Miller, rather than at tuples of constants; its types consist of closed conditions $\synt{\phi\geq r}$ in place of formulas; and its space is the space of Henkin types that satisfy the theory, with basic open sets given by finite theories, so that the hypothesis of Keisler's theorem becomes a statement about interiors in a \v{C}ech-complete space.

In the terms of Keisler and Miller, the hypothesis of Theorem~\ref{thm:B} says that no type of the sequence is locally realized.
They call a type locally realized when it is implied, through some tuple of terms, by a formula consistent with the theory~\cite[Definition~4.1]{Keisler-Miller:2001}; they show that every type realized in every model of a consistent theory is locally realized~\cite[Theorem~4.2]{Keisler-Miller:2001}, and they state local realization as the existence of a term mapping between spaces of types under which the preimage of the type has nonempty interior~\cite[Section~4]{Keisler-Miller:2001}.
Caicedo and Iovino use tuples of terms in the same way in their $\omega$-principal types~\cite[Definition~3.5]{Caicedo-Iovino:2014}, and Eagle does so in his principal types, crediting this use of terms to Keisler and Miller and to Caicedo and Iovino~\cite[Definition~4.10]{Eagle:2014}.

In Part~II we prove Scott's theorem for general structures (Theorem~\ref{thm:C}): every structure with countable sorts is characterized, among all structures, by a single sentence up to approximate back-and-forth equivalence.
This equivalence coincides with agreement on the full type, the values of all sentences of $\Loo$.
We give Scott's theorem a topological form (Theorem~\ref{thm:D}) and establish a version of López-Escobar's theorem that marks the reach of this form.
These theorems hold in the same generality as those of Part~I.

Scott proved that every countable structure for a countable vocabulary is characterized up to isomorphism, among countable structures, by a single sentence of $\Loo$~\cite{Scott:1965}.
The theorem has a topological reading: the structures whose universe is the set of natural numbers form a Polish space, the isomorphism classes are the orbits of the group of permutations of the natural numbers, and every orbit, as the set of models of a Scott sentence, is Borel.
López-Escobar proved the converse, that every Borel set of such structures invariant under the group is the set of models of a sentence of $\Loo$~\cite{LopezEscobar:1965}, and Vaught's transforms later gave the theorem a topological proof~\cite{Vaught:1974}.
Both theorems rest on equality: isomorphisms are bijections, and the group of permutations is the group of bijections of the universe.
Ben Yaacov, Doucha, Nies and Tsankov extended Scott's theorem, together with a version of López-Escobar's theorem, to separable metric structures whose predicates are uniformly continuous, with moduli prescribed by the signature, with respect to a complete metric~\cite{BenYaacov-Doucha-Nies-Tsankov:2017}; there the metric takes the part of equality, and structures without a metric lie outside the setting of their analysis.
General structures require neither equality nor a metric, and their sentences register a structure only up to approximation (Example~\ref{ex:dense}).
In Theorem~\ref{thm:C}, therefore, the place of isomorphism is taken by \emph{approximate equivalence}: two structures are approximately equivalent if, for every real number $r>0$, there is a back-and-forth system between them that matches the values of the atomic formulas to within $r$ times their a priori bounds.

\begin{maintheorem}
\label{thm:C}
Let a countable vocabulary and a structure $\fM$ whose sorts are countable be given.
There is a sentence $\sigma_\fM$ of $\Loo$ such that, for every structure $\fN$, the following conditions are equivalent: $\fN$ satisfies the closed condition $\synt{\sigma_\fM\leq0}$; the structures $\fM$ and $\fN$ are approximately equivalent; and the structures $\fM$ and $\fN$ have the same full type.
If, moreover, the vocabulary has a single sort and a predicate symbol for equality, and $\fM$ and $\fN$ are structures with countable sorts in which every atomic formula $\psi$ other than $\unity$ takes only the two values $r_\psi$ and $-r_\psi$, the predicate symbol for equality taking the value $1$ exactly on the diagonal, then $\fN$ satisfies $\synt{\sigma_\fM\leq0}$ exactly when $\fN$ is isomorphic to $\fM$.
\end{maintheorem}

The first assertion is Theorem~\ref{thm:scott}, and the second is Corollary~\ref{cor:scott-classical}.
The key observation is that the Scott formulas measure, by their values, the error with which a tuple of another structure matches a tuple of $\fM$ in the back-and-forth game, and that a strict inequality keeps this error from accumulating along the game.
Approximate equivalence is intrinsic to this syntax.
It can identify nonisomorphic countable structures even when both structures are complete, indeed discrete, for the pseudometrics induced by the formulas.
Thus passing to these completions does not turn the equivalence characterized by the Scott sentence into isomorphism (Example~\ref{ex:complete}).

Theorem~\ref{thm:C} has a topological form, in which the logic enters only the proof.
For a vocabulary with a single sort, let $\Xs$ be the Polish space of structures whose universe is the set $\NN$ of natural numbers.
The atomic diagram of a structure in $\Xs$ records the values of the atomic formulas at the tuples of natural numbers; surjections of $\NN$ onto $\NN$ act on these diagrams by reindexing them.
For structures $\fM$ and $\fN$ in $\Xs$, we compare their diagrams after reindexing each by a surjection.
The number $\delta(\fM,\fN)$ is the infimum, over all pairs of surjections, of the resulting uniform distances, normalized by the a priori bounds.
Two structures are \emph{exactly equivalent} if there is a back-and-forth system between them along which the values of the atomic formulas agree.

\begin{maintheorem}
\label{thm:D}
Let a countable vocabulary with a single sort be given, and let $\Xs$ and $\delta$ be as above.
\begin{enumerate}[(a)]
\item Two structures in $\Xs$ are approximately equivalent if and only if $\delta$ vanishes at them, and exactly equivalent if and only if their reindexed atomic diagrams coincide for some pair of surjections.
\item The function $\delta$ is a pseudometric on $\Xs$, and every class of approximate equivalence is the zero set of a bounded Borel function on $\Xs$ that is Lipschitz for $\delta$; in particular, every class is a Borel set, and two structures in $\Xs$ are approximately equivalent if and only if no Borel function on $\Xs$ that is Lipschitz for $\delta$ separates them.
\item A bounded function on $\Xs$ is Borel and takes the same value at exactly equivalent structures if and only if it is given by a sentence of the logic obtained from $\Loo$ by adding, as connectives, the continuous functions of finitely many real variables.
\end{enumerate}
\end{maintheorem}

Part~(a) is Lemma~\ref{lem:surjections}, part~(b) is Theorem~\ref{thm:topological-scott}, and part~(c) is Theorem~\ref{thm:lopez-escobar}.
The monoid of surjections, rather than the group of permutations, is what the absence of equality calls for: a structure is exactly equivalent to the structure obtained by giving each of its elements a twin, and no pair of permutations brings the two atomic diagrams within distance less than $1$ (Example~\ref{ex:twins}).
The main tool for part~(b) is Lemma~\ref{lem:transfer}, by which every sentence of $\Loo$ gives a function that is Lipschitz for $\delta$, and the main tool for part~(c) is the method of Vaught's transforms, carried out over the monoid of surjections.

Part~(c) marks the reach of part~(b).
A sentence of the enlarged logic can distinguish the value $\frac12$ of a predicate from values arbitrarily close to it, so some of the Borel functions invariant under exact equivalence separate approximately equivalent structures (Example~\ref{ex:exact}); replacing the functions of part~(b) by all such functions therefore changes the equivalence that they detect, from approximate equivalence to exact equivalence.
When the values of the atomic formulas range over finite sets separated by a uniform gap, the two equivalences coincide, and the invariant Borel functions are given by sentences of $\Loo$ itself (Theorem~\ref{thm:finite-values}); for two-valued structures, this is López-Escobar's theorem for classical $\Loo$ without equality, with the monoid of surjections in place of the group of permutations.

\subsection*{Organization of the paper}
Part~I consists of Sections~\ref{sec:topological-prelim}--\ref{sec:omit-type}.
In Section~\ref{sec:topological-prelim} we recall \v{C}ech's form of the Baire category theorem, and in Sections~\ref{sec:vocab-structures}--\ref{sec:type} we set up the structures, languages, fragments, theories and formal types.
In Section~\ref{sec:Henkin-type} we define Henkin types and prove Theorem~\ref{thm:A}.
In Section~\ref{sec:Henkin-topology} we show that, in the countable case, the space of Henkin types of a fragment is a $G_\delta$ subset of the product of intervals in which it lies, and Theorem~\ref{thm:B}, proved in Section~\ref{sec:omit-type}, then follows from \v{C}ech's form of the Baire category theorem.
Part~II consists of Sections~\ref{sec:scott} and~\ref{sec:topological-scott}.
In Section~\ref{sec:scott} we construct the Scott sentences and prove Theorem~\ref{thm:C}, and in Section~\ref{sec:topological-scott} we give the topological form of Theorem~\ref{thm:C} and prove Theorem~\ref{thm:D}.

\part{Model Existence and the Omitting Types Theorem}

\section{Topological Preliminaries}
	\label{sec:topological-prelim}

	The topological fact on which the paper rests is \v{C}ech's form of the Baire category theorem, which applies, in particular, to every countable intersection of open subsets of a compact Hausdorff space.
The theorem is applied in the space of Henkin types of a fragment (Section~\ref{sec:Henkin-type}), which, as we show in Section~\ref{sec:Henkin-topology}, is such an intersection when the vocabulary and the fragment are countable.
We recall the definitions that the theorem involves, with Engelking's \emph{General Topology}~\cite{Engelking:1989} as our reference for topology.

A topological space is a \textbf{Baire space} if the intersection of every sequence of dense open subsets of the space is dense; Baire spaces were introduced by Bourbaki in 1948 (see \cite[notes to Section~3.9]{Engelking:1989}).
The spaces to which this property will be applied are Tychonov spaces, that is, completely regular $T_1$-spaces; every subspace of a compact Hausdorff space is one (see \cite[Theorem~3.2.6]{Engelking:1989}).
Let $X$ be a Tychonov space.
A sequence $(\mathscr{U}_n)_{n<\omega}$ of open covers of $X$ is \textbf{complete} if every family $(F_i)_{i\in I}$ of closed subsets of $X$ that has the finite intersection property, and that contains, for every $n<\omega$, a member included in some element of $\mathscr{U}_n$, has nonempty intersection.
The space $X$ is \textbf{\v{C}ech-complete} if it has a complete sequence of open covers.
\v{C}ech~\cite{Cech:1937} introduced these spaces in 1937, under the name of topologically complete spaces, as the Tychonov spaces that are $G_\delta$ subsets of their Stone--\v{C}ech compactifications.
The characterization of \v{C}ech-complete spaces by sequences of open covers, adopted here as the definition, was established independently by Frol\'{\i}k in 1960 and by Arhangel'ski\u{\i} in 1961 (see \cite[Theorem~3.9.2 and the notes to Section~3.9]{Engelking:1989}).
By Theorems~3.9.1 and~3.9.2 of Engelking~\cite{Engelking:1989}, a Tychonov space is \v{C}ech-complete if and only if it is a $G_\delta$ subset of some, or equivalently of every, compactification of it.
Consequently, every $G_\delta$ subset $X$ of a compact Hausdorff space $K$ is \v{C}ech-complete.
The ambient compact spaces used below are products $\prod_{i\in I}[a_i,b_i]$ of real intervals with $a_i\le b_i$, endowed with the product topology; they are compact Hausdorff by Tychonov's theorem (see \cite[Theorem~3.2.4]{Engelking:1989}).
\v{C}ech also proved the following form of the Baire category theorem (see \cite[Corollary~3.9.4 and the notes to Section~3.9]{Engelking:1989}).

\begin{theorem}[Baire category theorem; \v{C}ech~\cite{Cech:1937}]
\label{thm:baire-category}
Every \v{C}ech-complete space is a Baire space.
\end{theorem}

	\section{Vocabularies and Structures}
\label{sec:vocab-structures}

The models whose existence the omitting types theorem asserts are real-valued structures.
Such a structure consists of one or more nonempty sets, called its sorts, of functions from finite products of sorts to sorts, and of its predicates.
The predicates of the structure are bounded real-valued functions on finite products of sorts.
Classical first-order structures for vocabularies without equality correspond to the case in which every predicate takes only the values $1$ and $-1$, read as true and false, respectively.

	Formulas speak about structures through a vocabulary, which is a set of symbols for sorts, functions and predicates.
A structure for the vocabulary interprets each symbol as an actual sort, function or predicate.
Since symbols must be kept apart from their interpretations, the symbols of a vocabulary are printed in \symb{this color}, and mathematical objects in black.
For the same reason, the variables, the connectives and the formal expressions built from them are printed in \synt{this color}.

	A \textbf{sort descriptor} is a purely syntactic symbol $\SM$.
Let $\MM$ be a nonempty set of sort descriptors.
A \textbf{function symbol for $\MM$} is a purely syntactic symbol $\Sf$ given together with
\begin{enumerate}[(i)]
\item a natural number $n=n_{\Sf}$, called the \textbf{arity of $\Sf$}, and
\item an $(n+1)$-tuple $(\SM_i)_{i\leq n}$ of sort descriptors in $\MM$, called the \textbf{extended arity of $\Sf$}.
\end{enumerate}
We indicate the arity and the extended arity of $\Sf$ by writing
\[
\Sf \colon \prod_{i<n} \SM_i \to \SM_n.
\]
A function symbol of arity $0$ is called a \textbf{constant symbol of sort $\SM_0$}.
A \textbf{predicate symbol for $\MM$} is a purely syntactic symbol $\Sp$ given together with
\begin{enumerate}[(i)]
\item a natural number $n=n_{\Sp}$, called the \textbf{arity of $\Sp$},
\item an $n$-tuple $(\SM_i)_{i< n}$ of sort descriptors in $\MM$, called the \textbf{extended arity of $\Sp$}, and
\item a real number $r = r_{\Sp} \ge 0$, called the \textbf{a priori bound of $\Sp$}.
\end{enumerate}
We indicate the arity, the extended arity and the a priori bound of $\Sp$ by writing
\[
\Sp \colon \prod_{i<n} \SM_i \to [-r,r].
\]
A \textbf{vocabulary} is a triple $V = (\MM, F, P)$ consisting of a nonempty set $\MM$ of sort descriptors, a set $F$ of function symbols for $\MM$, and a set $P$ of predicate symbols for $\MM$.
The set $P$ is required to contain a nullary predicate symbol $\unity$, called the \textbf{logical unity constant}, with a priori bound $r_{\unity} = 1$.
The vocabulary $V$ is \textbf{countable} if the sets $\MM$, $F$ and $P$ are countable.
The a priori bounds provide the compactness on which the topology of the paper rests.
Indeed, as we show in Section~\ref{sec:lang}, the value of a formula in a structure lies in an interval $[-r,r]$ determined by the formula alone.
Consequently, the types of all structures lie in a product of compact intervals, the compact Hausdorff space of formal types of Section~\ref{sec:type}.

	Given a vocabulary $V = (\MM, F, P)$, a \textbf{$V$-structure}, or real-valued structure with vocabulary $V$, is a triple
\[
\fM = \langle ({\SM}^\fM)_{{\SM}\in \MM}, ({\Sf}^\fM)_{{\Sf}\in F}, ({\Sp}^\fM)_{{\Sp}\in P} \rangle
\]
in which
\begin{itemize}
\item the interpretation ${\SM}^\fM$ of each sort descriptor $\SM\in \MM$ is a nonempty set,
\item the interpretation of each function symbol $\Sf\colon\prod_{i<n}\SM_i\to\SM_n$ in $F$ is a function
\[
{\Sf}^\fM \colon \prod_{i< n} {\SM}_i^\fM \to {\SM}_{n}^\fM,
\]
\item the interpretation of each predicate symbol $\Sp\colon\prod_{i<n}\SM_i\to[-r_{\Sp},r_{\Sp}]$ in $P$ is a function
\[
{\Sp}^\fM \colon \prod_{i< n} {\SM}_i^\fM \to [-r_{\Sp}, r_{\Sp}],
\]
and
\item the interpretation of the logical unity constant is ${\unity}^\fM = 1$.
\end{itemize}
Here a nullary symbol is interpreted on the one-point set $\prod_{i<0}\SM_i^\fM$, and its interpretation is identified with its value.

	\section{First-Order and Infinitary Languages}
\label{sec:lang}

Types, and later Henkin types, are formed from the values that structures assign to formulas.
We therefore build, from the symbols of a vocabulary, the terms and the formulas of a first-order language and of its infinitary extension, and interpret them in structures.

For the rest of the paper, we fix a vocabulary $V = (\MM, F, P)$.

	For each sort descriptor $\SM\in \MM$, we fix a sequence $({\Sx}_i^{\SM})_{i<\omega}$ of distinct symbols, not in $V$, chosen so that the sequences for distinct sort descriptors have no symbol in common.
Each symbol ${\Sx}_i^{\SM}$ is called a \textbf{variable of sort $\SM$}.

	\textbf{Terms of sort $\SM$} are defined recursively as follows:
\begin{itemize}
\item every variable of sort $\SM$ is a term of sort $\SM$;
\item for every function symbol $\Sf \colon \prod_{i<n}\SM_i\to \SM$ in $F$ and every tuple $\bar{\tau} = (\tau_i)_{i<n}$ of terms such that $\tau_i$ is of sort $\SM_i$ for each $i<n$, the expression $\Sf(\bar{\tau})$ is a term of sort $\SM$.
\end{itemize}
In particular, every constant symbol of sort $\SM$ is a term of sort $\SM$.
For a natural number $n$ and a tuple $(\SM_i)_{i<n}$ of sort descriptors in $\MM$, a tuple $(\tau_i)_{i<n}$ of terms such that $\tau_i$ is of sort $\SM_i$ for each $i<n$ is called a \textbf{term of sort $\prod_{i<n}\SM_i$}; for $n=0$, the only such term is the empty tuple.
A term is \textbf{closed} if no variable occurs in it.

	For a predicate symbol $\Sp \colon  \prod_{i<n}\SM_i \to [-r_{\Sp}, r_{\Sp}]$ in $P$ and a term $\bar{\tau}$ of sort $ \prod_{i<n}\SM_i$, the expression $\Sp(\bar{\tau})$ is called an \textbf{atomic $V$-formula} with \textbf{a priori bound} $r_{\Sp(\bar{\tau})} \coloneq r_{\Sp}$.
The other $V$-formulas are formed from the atomic ones by applying symbols called \textbf{logical connectives}:
\begin{itemize}
\item for each rational number $s$, the symbol $\Sneg$ is a logical connective;
\item the symbols $\Splus$ and $\Sor$ are logical connectives;
\item for each variable $\Sx$, the symbol $\supx$ is a logical connective.
\end{itemize}
The $V$-formulas and their a priori bounds are defined recursively, starting from the clause that every atomic $V$-formula is a $V$-formula.
If $\phi$ and $\psi$ are $V$-formulas with a priori bounds $r_\phi$ and $r_\psi$, respectively, then the expressions $\Sneg \phi$, $(\phi\Splus\psi)$, $(\phi\Sor\psi)$ and $\supx\phi$ are $V$-formulas for every rational number $s$ and every variable $\Sx$.
The a priori bounds of these four formulas are defined by
\[
\begin{gathered}
r_{\Sneg\phi} \coloneq |s| r_\phi, \qquad r_{\phi\Splus\psi} \coloneq r_\phi+r_\psi,\\
r_{\phi\Sor\psi} \coloneq \max\{r_\phi, r_\psi\}, \qquad r_{\supx\phi} \coloneq r_\phi,
\end{gathered}
\]
respectively.
Since the connectives $\Sneg$ and $\supx$ are written before the formula to which they apply, the parentheses in $(\phi\Splus\psi)$ and $(\phi\Sor\psi)$ ensure that every $V$-formula is formed in only one way.
Hence the a priori bound of a $V$-formula, like the other notions defined below by recursion on formulas, is well defined; we omit the parentheses when no ambiguity arises.
The set of all $V$-formulas is called the \textbf{first-order language for $V$} and is denoted by $L^V$.

	The set $\free(\xi)$ of \textbf{variables free in a $V$-formula $\xi$} is defined recursively as follows:
\begin{itemize}
\item if $\xi$ is atomic, of the form $\Sp(\bar{\tau})$, then $\free(\xi)$ is the set of variables that occur in the terms of the tuple $\bar{\tau}$;
\item if $\xi$ is of the form $\Sneg \phi$, then $\free(\xi) \coloneq \free(\phi)$;
\item if $\xi$ is of the form $\phi\Splus\psi$ or $\phi\Sor\psi$, then $\free(\xi) \coloneq \free(\phi)\cup\free(\psi)$;
\item if $\xi$ is of the form $\supx \phi$, then $\free(\xi) \coloneq \free(\phi)\setminus\{\Sx\}$.
\end{itemize}
A \textbf{$V$-sentence} is a $V$-formula with no free variables, that is, a $V$-formula $\phi$ such that $\free(\phi) = \emptyset$.

	For a term $\tau$ and a tuple $\bar{x}$ of distinct variables, we write $\tau(\bar{x})$ to indicate that every variable occurring in $\tau$ is among $\bar{x}$.
Similarly, for a formula $\phi$, we write $\phi(\bar{x})$ to indicate that every variable free in $\phi$ is among $\bar{x}$.
If $\bar{x} = (\Sx_i)_{i<n}$, where $\Sx_i$ is a variable of sort $\SM_i$ for each $i<n$, we write
\[
{}_{\bar{x}} \SM \coloneq \prod_{i<n}\SM_i
\]
and say that $\bar{x}$ is a variable of sort ${}_{\bar{x}}\SM$.
For a $V$-structure $\fM$, we write
\[
{}_{\bar{x}} \SM^{\fM} \coloneq \prod_{i<n}\SM_i^\fM.
\]

	The infinitary language $\Loo^V$ is obtained by adding to the formation rules of $L^V$ a new logical connective $\join$, called the \textbf{infinitary connective of countable join}.
More precisely, the language $\Loo^V$, together with the a priori bounds and the free variables of its members, is defined by recursion on the countable ordinals.
Let $\cL_0$ be the set of atomic $V$-formulas, and let $\cL_\lambda \coloneq \bigcup_{\alpha<\lambda}\cL_\alpha$ for every countable limit ordinal $\lambda$.
For $\alpha<\omega_1$, the set $\cL_{\alpha+1}$ is obtained by adding to $\cL_\alpha$ two kinds of expressions formed from its members.
The first are the expressions $\Sneg\phi$, $\phi\Splus\psi$, $\phi\Sor\psi$ and $\supx\phi$ with $\phi,\psi\in\cL_\alpha$, for all rational numbers $s$ and all variables $\Sx$; their a priori bounds and free variables are given by the clauses stated for $L^V$.
The second concern the sequences $\Phi = (\phi_n)_{n<\omega}$ of members of $\cL_\alpha$ for which $\sup_{n<\omega} r_{\phi_n}$ is finite and, for some finite tuple $\bar{x}$ of variables, every $\phi_n$ is of the form $\phi_n(\bar{x})$.
For every such sequence, the expression $\join\Phi$ is added, with a priori bound $r_{\join\Phi} \coloneq \sup_{n<\omega} r_{\phi_n}$.
The free variables of $\join\Phi$ are those of the formulas $\phi_n$, that is, $\free(\join\Phi) \coloneq \bigcup_{n<\omega}\free(\phi_n)$, a subset of the finite set of variables of $\bar{x}$.
Finally, we set $\Loo^V \coloneq \bigcup_{\alpha<\omega_1}\cL_\alpha$.
Since $\omega_1$ is regular, every sequence $(\phi_n)_{n<\omega}$ of members of $\Loo^V$ consists of members of some $\cL_\alpha$, and the set $\Loo^V$ is therefore closed under all these formation rules.

We also write $\join_n \phi_n$ for $\join\Phi$, and we denote the set of sentences of $\Loo^V$ by $\Soo^V$.
A formula of $\Loo^V$ may have infinite length, but only finitely many free variables.
The two languages differ in size: if $V$ is countable, the first-order language $L^V$ is countable, whereas the infinitary language $\Loo^V$ has the cardinality of the continuum.
For this reason, the Baire category argument of this paper takes place in countable fragments (see Section~\ref{sec:fragments}).

	Let $\fM$ be a $V$-structure, and let $\tau(\bar{x})$ be a term of sort $\SM$, where $\bar{x} = (\Sx_j)_{j<m}$.
The \textbf{interpretation of $\tau(\bar{x})$ in $\fM$} is the function $\tau^\fM \colon {}_{\bar{x}}\SM^\fM \to \SM^\fM$ defined recursively, for $\bar{a} = (a_j)_{j<m}$ in ${}_{\bar{x}}\SM^\fM$, as follows:
\begin{itemize}
\item the interpretation of the variable $\Sx_j$ is the projection $\Sx_j^\fM(\bar{a}) \coloneq a_j$;
\item if $\Sf \colon \prod_{i<n}\SM_i \to \SM$ is a function symbol and $\bar{\tau} = (\tau_i)_{i<n}$ is a term of sort $\prod_{i<n}\SM_i$, then the interpretation of the term $\Sf(\bar{\tau})$ is given by $(\Sf(\bar{\tau}))^\fM(\bar{a}) \coloneq \Sf^\fM ((\tau_i^\fM(\bar{a}))_{i<n})$.
\end{itemize}
The interpretation of a closed term $\tau$ is a function on the one-point set ${}_{\bar{x}}\SM^\fM$ with $\bar{x}$ empty, and is identified with its value, an element of $\SM^\fM$.
The \textbf{interpretation of a formula $\xi(\bar{x})$ of $\Loo^V$ in $\fM$} is the function $\xi^\fM$ on ${}_{\bar{x}}\SM^\fM$ defined recursively, for $\bar{a}$ in ${}_{\bar{x}}\SM^\fM$, as follows:
\begin{itemize}
\item if $\xi$ is atomic, of the form $\Sp(\bar{\tau})$ with $\bar{\tau} = (\tau_i)_{i<n}$, then
\[
\xi^\fM(\bar{a}) \coloneq \Sp^\fM ((\tau_i^\fM(\bar{a}))_{i<n});
\]
\item if $\xi$ is of the form $\Sneg \phi$, where $s$ is a rational number and $\phi\in \Loo^V$, then
\[
\xi^\fM(\bar{a}) \coloneq s\cdot \phi^\fM(\bar{a});
\]
\item if $\xi$ is of the form $\phi \Splus \psi$, where $\phi,\psi\in\Loo^V$, then
\[
\xi^\fM(\bar{a}) \coloneq \phi^\fM(\bar{a})+\psi^\fM(\bar{a});
\]
\item if $\xi$ is of the form $\phi \Sor \psi$, where $\phi, \psi\in \Loo^V$, then
\[
\xi^\fM(\bar{a}) \coloneq \max\{\phi^\fM(\bar{a}), \psi^\fM(\bar{a})\};
\]
\item if $\xi$ is of the form $\supy \phi$, where $\phi\in \Loo^V$, then
\[
\xi^\fM(\bar{a}) \coloneq \sup_{b\in {}_{\Sy} \SM^\fM} \phi^\fM(\bar{a}[b/\Sy]),
\]
where, if $\Sy$ is not among $\bar{x}$, the formula $\phi$ is regarded as a formula $\phi(\bar{x}\Sy)$, and $\bar{a}[b/\Sy]$ is the tuple $\bar{a}b$; if $\Sy$ is the variable $\Sx_j$ of $\bar{x}$, the formula $\phi$ is regarded as a formula $\phi(\bar{x})$, and $\bar{a}[b/\Sy]$ is the tuple obtained from $\bar{a}$ by replacing its entry $a_j$ by $b$;
\item if $\xi$ is of the form $\join\Phi$, where $\Phi = (\phi_n(\bar{x}))_{n<\omega}$ is a sequence of formulas of $\Loo^V$, then
\[
\xi^\fM(\bar{a}) \coloneq \sup_{n<\omega} \phi_n^\fM(\bar{a}).
\]
\end{itemize}
An induction on formulas shows that the suprema in this definition are finite and that $\xi^\fM$ takes its values in $[-r_\xi, r_\xi]$.
When $\xi$ is a sentence, the tuple $\bar{x}$ may be taken empty, and $\xi^\fM$ is identified with its value, a real number in $[-r_\xi, r_\xi]$.

\enlargethispage{2pt}%
The finite connectives suffice to approximate every continuous function of finitely many formula values on their prescribed intervals.
Indeed, for a positive integer $k$ and formulas $\phi_1,\dots,\phi_k$, write $K=\prod_{i=1}^k[-r_{\phi_i},r_{\phi_i}]$.
The coordinate functions on the compact box $K$, together with the constant function $1$, generate under addition, rational scalar multiplication and maximum a rational vector lattice of continuous functions that separates points.\footnote{Any fixed countable dense set of real scalar coefficients suffices for this approximation argument; the rational coefficients are a convenient choice. The logical unity constant already supplies every rational constant as $q\cdot\unity$, with $q\in\QQ$.}
Its uniform closure is a real vector lattice containing the constants, so the lattice Stone--Weierstrass theorem makes it all of $C(K,\RR)$, the space of continuous real-valued functions on $K$.
Thus, for every continuous function $f\colon K\to\RR$ and every real number $\varepsilon>0$, there is an expression $p$ in these finite connectives such that $|f(t)-p(t)|<\varepsilon$ for every $t\in K$.
In particular, the approximation is uniform over all structures when $t$ is the tuple of values of $\phi_1,\dots,\phi_k$.

Several proofs below proceed by induction on formulas of infinite length, and for these we use the \textbf{rank} of a formula, defined recursively as follows.
Atomic formulas have rank $0$; the rank of $\Sneg\phi$ and of $\supx\phi$ is the rank of $\phi$ plus one; the rank of $\phi\Splus\psi$ and of $\phi\Sor\psi$ is the larger of the ranks of $\phi$ and $\psi$ plus one; and the rank of $\join\Phi$ is the supremum of the ranks of the members of $\Phi$ plus one.
Every rank is a countable ordinal, and the rank of a proper subformula of a formula is smaller than the rank of the formula.
Since substituting terms for variables changes only the atomic subformulas, it preserves the rank.
These definitions apply to every vocabulary, and in particular to the expanded vocabulary of Section~\ref{sec:expand-const}.

	\section{Expansion by Constants}
\label{sec:expand-const}

In the method introduced by Henkin in 1949, the elements of a structure are named by new constant symbols; the first step is to add these symbols to the fixed vocabulary.
Let $C = \{ {\Sc}_n^{\SM} \colon n<\omega,\ \SM \in \MM \}$ be a set of distinct constant symbols, not in $V$, such that ${\Sc}_n^{\SM}$ is of sort $\SM$ for every $n<\omega$ and every $\SM\in\MM$.
The \textbf{expansion of $V$ by $C$} is the vocabulary $V_C \coloneq (\MM, F\cup C, P)$, which is countable if $V$ is; the infinitary language $\Loo^{V_C}$ is denoted by $\Loo^C$, and its set of sentences by $\Soo^C$.

	Given a $V$-structure $\fM = \langle (\SM^\fM)_{\SM\in\MM}, (\Sf^\fM)_{\Sf\in F}, (\Sp^\fM)_{\Sp\in P} \rangle$, an \textbf{expansion of $\fM$ to $V_C$} is a $V_C$-structure
\[
\fN = \langle (\SM^\fN)_{\SM\in\MM}, (\Sf^\fN)_{\Sf\in F\cup C}, (\Sp^\fN)_{\Sp\in P} \rangle
\]
that agrees with $\fM$ on the symbols of $V$, in the sense that $\SM^\fN = \SM^\fM$ for every $\SM\in\MM$, that $\Sf^\fN = \Sf^\fM$ for every $\Sf\in F$, and that $\Sp^\fN = \Sp^\fM$ for every $\Sp\in P$.
Each expansion $\fN$ of $\fM$ to $V_C$ is determined by the elements $({\Sc}_n^{\SM})^\fN$ of $\SM^\fM$ that interpret the new constant symbols.
Conversely, every $V_C$-structure $\fN$ is an expansion of exactly one $V$-structure, its \textbf{$V$-reduct}, obtained by forgetting the interpretations of the symbols in $C$.

	In the structures built by the Henkin construction, every element has a name.
A $V_C$-structure $\fN$ is \textbf{canonical} if, for every sort descriptor $\SM\in\MM$, every element of $\SM^\fN$ is the interpretation in $\fN$ of a closed $V_C$-term of sort $\SM$.

In a canonical structure, a supremum over a sort is a supremum over the values of closed terms, and it can be computed from sentences by substituting these terms for the variable.
The following lemma, which justifies this computation, holds in every $V_C$-structure.
For a formula $\psi(\bar{z}\bar{x})$ and a tuple $\bar{\tau}$ of closed terms, we write $\psi(\bar{z}\bar{\tau})$ for the formula obtained by substituting the terms of $\bar{\tau}$ for the free occurrences of the corresponding variables of $\bar{x}$; a closed term is free for every variable.

\begin{lemma}[Substitution]
	\label{lem:substitution}
	Let $\fN$ be a $V_C$-structure, let $\psi(\bar{z}\bar{x})$ be a formula of $\Loo^C$, where $\bar{x}$ is a tuple of distinct variables not among $\bar{z}$, and let $\bar{\tau}$ be a tuple of closed $V_C$-terms of sort ${}_{\bar{x}}\SM$.
	Then $\psi(\bar{z}\bar{\tau})^\fN(\bar{b}) = \psi^\fN(\bar{b}\,\bar{\tau}^\fN)$ for every tuple $\bar{b}\in{}_{\bar{z}}\SM^\fN$.
\end{lemma}

\begin{proof}
	The corresponding identity for terms is proved by induction on terms, and the lemma follows by induction on the rank of $\psi$.
	The parameters $\bar{z}$ are needed in the step for suprema, in which the variable of the supremum joins them.
\end{proof}

	\section{Fragments}
\label{sec:fragments}

Countability is indispensable to the Baire category argument, which writes the Henkin types of Section~\ref{sec:Henkin-type} as a countable intersection of open sets attached to sentences (see Section~\ref{sec:Henkin-topology}).
The infinitary language, however, has the cardinality of the continuum even when the vocabulary is countable (see Section~\ref{sec:lang}).
The argument is therefore carried out in fragments: sets of formulas that contain the first-order language and are closed under the first-order connectives, under subformulas and under substitution of terms.
A fragment generated by countably many formulas over a countable vocabulary is countable, as we show below.

Let $\xi$ be a formula of $\Loo^V$.
The set $\sub(\xi)$ of \textbf{subformulas of $\xi$} is defined recursively as follows:
\begin{itemize}
\item if $\xi$ is atomic, we set $\sub(\xi) \coloneq \{\xi\}$;
\item if $\xi$ is of the form $\Sneg \phi$, we set $\sub(\xi) \coloneq \{\xi\} \cup \sub(\phi)$;
\item if $\xi$ is of the form $\phi\Splus\psi$ or $\phi \Sor \psi$, we set $\sub(\xi) \coloneq \{\xi\} \cup \sub(\phi) \cup \sub(\psi)$;
\item if $\xi$ is of the form $\supx \phi$, we set $\sub(\xi) \coloneq \{\xi\} \cup \sub(\phi)$;
\item if $\xi$ is of the form $\join_n \phi_n$, we set $\sub(\xi) \coloneq \{\xi\} \cup \bigcup_{n<\omega}\sub(\phi_n)$.
\end{itemize}
	A \textbf{fragment} of $\Loo^V$ is a set $\frag \subseteq \Loo^V$ with the following properties:
\begin{enumerate}[(F1)]
\item it contains $L^V$;
\item it is closed under subformulas, that is, $\sub(\phi)\subseteq \frag$ for every $\phi\in\frag$;
\item it is closed under the first-order connectives, that is, if $\phi, \psi\in \frag$, then $\Sneg \phi \in \frag$ for every rational number $s$, $\phi\Splus\psi\in\frag$, $\phi \Sor\psi\in \frag$, and $\supx \phi\in\frag$ for every variable $\Sx$;
\item it is closed under substitution of terms, that is, if $\phi(\bar{x}\Sy)\in\frag$ and $\tau$ is a term of sort $\sorty$ that is free for $\Sy$ in $\phi$ (no free occurrence of $\Sy$ in $\phi$ lies within a subformula $\synt{\sup_z}\psi$ of $\phi$ such that $\synt{z}$ occurs in $\tau$), then the formula $\phi(\bar{x}\tau)$, obtained from $\phi$ by substituting $\tau$ for the free occurrences of $\Sy$, belongs to $\frag$.
\end{enumerate}
	The intersection of a nonempty family of fragments of $\Loo^V$ is a fragment, and the language $\Loo^V$ is itself a fragment.
Hence, for every set $\Phi\subseteq\Loo^V$, the intersection of all fragments containing $\Phi$ is the smallest fragment containing $\Phi$, which is called the \textbf{fragment generated by $\Phi$}.
If $V$ and $\Phi$ are countable, then the fragment generated by $\Phi$ is countable.
Indeed, let $\frag_0 \coloneq L^V\cup\bigcup_{\phi\in\Phi}\sub(\phi)$; for each $k<\omega$, let $\frag_{k+1}$ be obtained by adding to $\frag_k$ the formulas that (F2)--(F4) require for members of $\frag_k$.
Since every formula has countably many subformulas, and since there are countably many rational numbers, variables and terms, an induction on $k$ shows that every set $\frag_k$ is countable.
The union of the increasing sequence $(\frag_k)_{k<\omega}$ is a fragment that contains $\Phi$ and is contained in every fragment containing $\Phi$.
Consequently, the fragment generated by $\Phi$ is this union, which is countable.

	Let $\frag$ be a fragment of $\Loo^V$.
Since $\Loo^V$ is contained in the infinitary language $\Loo^C$ of the expanded vocabulary $V_C$ of Section~\ref{sec:expand-const}, we may consider the fragment of $\Loo^C$ generated by $\frag$, which we denote by $\frag_C$.
If $V$ and $\frag$ are countable, then so is $\frag_C$, by the preceding paragraph applied to the countable vocabulary $V_C$.
We denote the set of sentences of $\frag$ by $\SF$ and the set of sentences of $\frag_C$ by $\SFC$.

	\section{Conditions and Theories}
\label{sec:cond-theories}

In real-valued logic, the role of the statements of classical logic is played by conditions, formal strict inequalities between a sentence and a rational number.
Sets of conditions, called theories, specify the models to be constructed, whereas the types to be omitted are expressed by closed conditions.
Strict conditions are taken as the basic ones because each of them defines an open set in the topology of Section~\ref{sec:Henkin-topology}, as the Baire category argument requires.

A \textbf{lower-bound condition} is a formal expression $\synt{\phi > r}$, and an \textbf{upper-bound condition} is a formal expression $\synt{\phi < r}$.
Each of these expressions is formed from a sentence $\phi$ of the infinitary language $\Loo^C$ of the expanded vocabulary $V_C$ (see Section~\ref{sec:expand-const}) and a rational number $r$.
A \textbf{condition} is a lower-bound condition or an upper-bound condition, and a \textbf{theory} is a set of conditions.
For a set $\Sigma$ of sentences, a \textbf{$\Sigma$-theory} is a theory whose conditions are on sentences in $\Sigma$.

	For a rational number $r$, we write
\begin{itemize}
\item $\synt{\phi \geq r}$ for the countable theory $\{ \synt{\phi>r-\frac{1}{n}} \colon n\in \ZZ_+ \}$, and
\item $\synt{\phi \leq r}$ for the countable theory $\{ \synt{\phi<r+\frac{1}{n}} \colon n\in \ZZ_+ \}$.
\end{itemize}
The expressions $\synt{\phi\geq r}$ and $\synt{\phi\leq r}$ are called \textbf{closed conditions}.

	Let $\fN$ be a $V_C$-structure.
We say that $\fN$ \textbf{satisfies} the condition $\synt{\phi>r}$, and write $\fN \models \synt{\phi>r}$, if $\phi^\fN>r$; it satisfies the condition $\synt{\phi<r}$ if $\phi^\fN<r$.
A $V_C$-structure that satisfies every condition of a theory is called a \textbf{model} of the theory; the same term is used for $V$-structures and theories whose conditions are on $V$-sentences.
Thus a structure satisfies the closed condition $\synt{\phi\geq r}$ exactly when the value of $\phi$ in it is at least $r$.

	Under the reading of $1$ as true and $-1$ as false, the connectives generalize those of classical first-order logic, as the following equivalences show.
Each of these equivalences holds for every $V_C$-structure $\fN$ and every rational number $r$.
Multiplication by $-1$ plays the role of negation: for a sentence $\phi$, if $\synt{-\phi}$ denotes the sentence obtained by applying the connective $\Sneg$ with $s=-1$, then
\begin{align*}
\fN \models \synt{-\phi>r} \quad &\text{iff} \quad (-\phi)^\fN > r\\
&\text{iff} \quad -(\phi^\fN)>r \\
&\text{iff} \quad \phi^\fN < -r.
\end{align*}
In particular, for $r=0$, a $V_C$-structure satisfies $\synt{-\phi>0}$ exactly when it satisfies $\synt{\phi<0}$.
The binary connective $\Sor$ is the counterpart of disjunction: for sentences $\phi$ and $\psi$, we have
\begin{align*}
\fN \models \synt{\phi\Sor\psi > r} \quad &\text{iff} \quad \max\{\phi^\fN, \psi^\fN\} > r\\
&\text{iff} \quad \fN\models\synt{\phi>r} \text{ or } \fN\models\synt{\psi>r}.
\end{align*}
The supremum $\supx$ takes the place of the existential quantifier: for a formula $\phi(\Sx)$ of $\Loo^C$, the condition $\synt{\supx\phi>r}$ is characterized as follows:
\begin{align*}
\fN \models \synt{\supx\phi>r} \quad &\text{iff} \quad (\supx\phi)^\fN > r\\
&\text{iff} \quad \sup \left\{\phi^\fN(a) \colon a\in \sortx^\fN \right\} > r \\
&\text{iff} \quad \phi^\fN(a) > r \text{ for some } a\in \sortx^\fN.
\end{align*}
Finally, the join $\join$ generalizes countable disjunction: for a sequence $\Phi$ of sentences of $\Loo^C$ whose join is a sentence of $\Loo^C$, the condition $\synt{\join\Phi>r}$ is characterized in the same way:
\begin{align*}
\fN \models \synt{\join \Phi > r} \quad &\text{iff} \quad \left(\join \Phi \right)^\fN > r\\
&\text{iff} \quad \sup \left\{\phi^\fN \colon \phi\in \Phi \right\} > r \\
&\text{iff} \quad \phi^\fN > r \text{ for some } \phi\in \Phi.
\end{align*}

	\section{Types of Structures and Formal Types}
\label{sec:type}

The type of a structure records the values that the structure assigns to a set of sentences; a formal type is a family of values in the prescribed intervals, whether or not these values come from a structure.
The Henkin conditions of Section~\ref{sec:Henkin-type} select formal types that do come from structures.
These Henkin types form the subspace of the compact space of formal types in which the Baire category argument takes place.

Let $\Phi$ be a set of sentences of the infinitary language $\Loo^C$ of the expanded vocabulary $V_C$.
We write $\prod\Phi$ for the product $\prod_{\phi\in\Phi}[-r_\phi, r_\phi]$ endowed with the product topology, where $r_\phi$ is the a priori bound of $\phi$.
In this topology, the product $\prod\Phi$ is a compact Hausdorff space (see Section~\ref{sec:topological-prelim}); its elements are called \textbf{formal $\Phi$-types}.
For a formal $\Phi$-type $\fm$ and a sentence $\phi\in\Phi$, we write $\phi^\fm$ for the $\phi$-coordinate of $\fm$ and say that $\fm$ satisfies a condition $\synt{\phi>r}$ or $\synt{\phi<r}$ if $\phi^\fm>r$ or $\phi^\fm<r$, respectively.
If $\fN$ is a $V_C$-structure, the \textbf{$\Phi$-type of $\fN$} is the formal $\Phi$-type
\[
\Phitp(\fN) \coloneq (\phi^\fN)_{\phi\in\Phi}.
\]
When $\Phi$ is contained in the set $\Soo^V$ of sentences of $\Loo^V$, the same definition applies to a $V$-structure $\fM$.
The $\Soo^V$-type $\fulltp(\fM)$ is called the \textbf{full type} of $\fM$; similarly, the full type of a $V_C$-structure is its type for the set $\Soo^C$ of sentences of $\Loo^C$.

	\section{Henkin Types}
	\label{sec:Henkin-type}

In his proof of the completeness theorem, Henkin~\cite{Henkin:1949} built a model of a consistent theory out of new individual constants, chosen to provide every existential sentence of the theory with a witness among them.
The Henkin conditions below demand of a formal type that its values respect the connectives and, in the real-valued form of Henkin's requirement, that the value of a supremum over a variable be the supremum of the values obtained by substituting closed terms for that variable.
In this section we show that the formal types satisfying the Henkin conditions are exactly the types of canonical structures, structures in which every element is the value of a closed term.
With the help of a downward L\"owenheim--Skolem theorem for countable fragments, we then determine which formal types are realized by structures.

Let $\cG$ be a fragment of the infinitary language $\Loo^C$ of the expanded vocabulary $V_C$, and let $S^{\cG}$ be the set of its sentences.
Let $\Sx$ be a variable, let $\phi$ be a formula whose free variables are among $\Sx$, and let $\tau$ be a closed $V_C$-term of sort $\sortx$.
We write $\phi(\tau)$ for the sentence obtained by substituting $\tau$ for the free occurrences of $\Sx$ in $\phi$, as in condition (F4) of the definition of a fragment (see Section~\ref{sec:fragments}).
A formal type $\fm\in\prod S^{\cG}$ is an \textbf{$S^{\cG}$-Henkin type} if it satisfies the following \textbf{Henkin conditions}:
\begin{enumerate}[(H1)]
	\item $\unity^\fm = 1$;
	\item $(\Sneg \phi)^\fm = s\,\phi^\fm$ for every $\phi\in S^{\cG}$ and every $s\in \QQ$;
	\item $(\phi\Splus\psi)^\fm = \phi^\fm+\psi^\fm$ for all $\phi,\psi\in S^{\cG}$;
	\item $(\phi\Sor\psi)^\fm = \max\{\phi^\fm, \psi^\fm\}$ for all $\phi, \psi\in S^{\cG}$;
	\item $(\supx\phi)^\fm = \sup\{\phi(\tau)^\fm \colon \tau \text{ is a closed } V_C\text{-term of sort } \sortx \}$ for every variable $\Sx$ and every formula $\phi\in\cG$ whose free variables are among $\Sx$;
	\item $\left(\join \Phi\right)^\fm = \sup\{\phi^\fm \colon \phi\in \Phi\}$ for every sentence $\join \Phi \in S^{\cG}$.
\end{enumerate}
Every sentence that occurs in these conditions belongs to $S^{\cG}$: the logical unity constant $\unity$ by (F1), and the sentences $\Sneg\phi$, $\phi\Splus\psi$, $\phi\Sor\psi$ and $\supx\phi$ by (F3).
The members of $\Phi$ belong to $S^{\cG}$ by (F2), since the free variables of $\join\Phi$ are those of its members; the sentences $\phi(\tau)$ belong to $S^{\cG}$ by (F4).
The language $\Loo^C$ is a fragment of itself, whose set of sentences is $\Soo^C$; its $\Soo^C$-Henkin types lie in $\prod\Soo^C$ and are simply called \textbf{Henkin types}.

To pass from formal types back to structures, we follow Henkin and take as elements the closed terms themselves; the structure so obtained depends only on the values of the atomic sentences.
Let $\cA$ be the set of atomic $V_C$-sentences; by (F1), it is contained in $S^{\cG}$ for every fragment $\cG$ of $\Loo^C$.
An \textbf{atomic type} is a formal type $\fa\in\prod\cA$ such that $\unity^\fa = 1$.
Every atomic type $\fa$ determines a $V_C$-structure $\fN_\fa$, the \textbf{term structure of $\fa$}, as follows.
For each $\SM\in\MM$, the set $\SM^{\fN_\fa}$ consists of the closed $V_C$-terms of sort $\SM$; it is nonempty, since it contains the constants $\Sc_n^{\SM}$.
Each function symbol $\Sf$ of $V_C$ is interpreted syntactically, by $\Sf^{\fN_\fa}(\bar{\tau}) \coloneq \Sf(\bar{\tau})$, and each predicate symbol $\Sp$ by $\Sp^{\fN_\fa}(\bar{\tau}) \coloneq \Sp(\bar{\tau})^\fa$.
The interpretation $\Sp^{\fN_\fa}$ is a function into $[-r_{\Sp}, r_{\Sp}]$, because the a priori bound of the atomic sentence $\Sp(\bar{\tau})$ is $r_{\Sp}$.
Since, moreover, the logical unity constant satisfies $\unity^{\fN_\fa} = \unity^\fa = 1$, the term structure $\fN_\fa$ is a $V_C$-structure.
No quotient of the set of terms is needed, since the vocabulary has no equality symbol.
An induction on terms shows that $\tau^{\fN_\fa} = \tau$ for every closed $V_C$-term $\tau$; hence the term structure $\fN_\fa$ is canonical and satisfies $\cA\text{-}\!\tp(\fN_\fa) = \fa$.

\begin{lemma}
	\label{lem:henkin-unique}
	Let $\cG$ be a fragment of $\Loo^C$.
	\begin{enumerate}[(a)]
		\item If $\fN$ is a canonical $V_C$-structure, then $S^{\cG}\text{-}\!\tp(\fN)$ is an $S^{\cG}$-Henkin type.
		\item Two $S^{\cG}$-Henkin types that agree on $\cA$ are equal.
		\item A formal type $\fm\in\prod S^{\cG}$ is the $S^{\cG}$-type of a canonical $V_C$-structure if and only if it is an $S^{\cG}$-Henkin type; in that case, it is the $S^{\cG}$-type of the term structure of its restriction to $\cA$.
	\end{enumerate}
\end{lemma}

\begin{proof}
	The idea is that the Henkin conditions repeat the recursive definition of the values of sentences, except that a supremum over a sort becomes a supremum over closed terms.
	The Henkin conditions therefore determine the value of every sentence from the values of the atomic ones, which is part (b).
	Part (a) asserts that the same conditions hold in a canonical structure, in which every element is the value of a closed term and the two suprema coincide.

	(a) Conditions (H1), (H2), (H3), (H4) and (H6) hold for the values $\sigma^\fN$ of the sentences $\sigma\in S^{\cG}$, by the definition of a $V_C$-structure and the recursive definition of the interpretation of formulas.
	For (H5), let $\phi\in\cG$ be a formula whose free variables are among $\Sx$.
	By the substitution lemma (Lemma~\ref{lem:substitution}), applied with $\bar{z}$ empty, we have $\phi(\tau)^\fN = \phi^\fN(\tau^\fN)$ for every closed $V_C$-term $\tau$ of sort $\sortx$.
	Since $\fN$ is canonical, every element of $\sortx^\fN$ is of the form $\tau^\fN$ for a closed $V_C$-term $\tau$ of sort $\sortx$; this gives
	\[
	(\supx\phi)^\fN = \sup_{b\in\sortx^\fN}\phi^\fN(b) = \sup_{\tau}\phi^\fN(\tau^\fN) = \sup_{\tau}\phi(\tau)^\fN .
	\]

(b) Let $\fm$ and $\fm'$ be $S^{\cG}$-Henkin types that agree on $\cA$.
	We show by induction on the rank of $\sigma$ that $\sigma^\fm = \sigma^{\fm'}$ for every $\sigma\in S^{\cG}$, the atomic case being the hypothesis.
	Suppose that $\sigma$ is of the form $\Sneg\phi$, $\phi\Splus\psi$, $\phi\Sor\psi$ or $\join\Phi$.
	In that case, the formulas $\phi$, $\psi$ and the members of $\Phi$ are sentences in $S^{\cG}$ of smaller rank, on which the types $\fm$ and $\fm'$ agree by the induction hypothesis.
	Conditions (H2), (H3), (H4) and (H6), respectively, compute $\sigma^\fm$ and $\sigma^{\fm'}$ alike from the values of these sentences.
	Suppose, finally, that $\sigma$ is of the form $\supx\phi$.
	For every closed $V_C$-term $\tau$ of sort $\sortx$, the sentence $\phi(\tau)$ then belongs to $S^{\cG}$ and has smaller rank.
	The types $\fm$ and $\fm'$ agree on each of these sentences for the same reason, and condition (H5) therefore gives $\sigma^\fm = \sigma^{\fm'}$.

	(c) By (a), the $S^{\cG}$-type of a canonical $V_C$-structure is an $S^{\cG}$-Henkin type.
	Conversely, let $\fm$ be an $S^{\cG}$-Henkin type, and let $\fa$ be its restriction to $\cA$, an atomic type by (H1).
	The term structure $\fN_\fa$ is canonical, with $\cA\text{-}\!\tp(\fN_\fa) = \fa$.
	By (a), the $S^{\cG}$-type of $\fN_\fa$ is an $S^{\cG}$-Henkin type that agrees with $\fm$ on $\cA$; by (b), the two types are equal.
\end{proof}

Let $\fa$ be an atomic type.
By Lemma~\ref{lem:henkin-unique}(a), applied to $\Loo^C$, the full type $\Soo^C\text{-}\!\tp(\fN_\fa)$ is a Henkin type that extends $\fa$.
Part (b) of that lemma shows that no other Henkin type extends $\fa$.
We call this full type the \textbf{Henkin type of $\fa$} and denote it by $\hat{\fa}$.
For every fragment $\cG$ of $\Loo^C$, the same lemma shows that the restriction of $\hat{\fa}$ to $S^{\cG}$ is the only $S^{\cG}$-Henkin type that extends $\fa$.

Let $\fN$ be a $V_C$-structure, and let $\fa \coloneq \cA\text{-}\!\tp(\fN)$.
If $\fN$ is canonical, then Lemma~\ref{lem:henkin-unique}(a) shows that the full type $\Soo^C\text{-}\!\tp(\fN)$ is a Henkin type that extends $\fa$, and therefore equals $\hat{\fa}$.
In other words, we have $\phi^\fN = \phi^{\hat{\fa}}$ for every $\phi \in \Soo^C$.
If $\fN$ is not canonical, the two values may differ, and in either direction.
For example, suppose that $V$ has a single sort descriptor $\SM$ and no function symbols, and that the only predicate symbol of $V$ besides $\unity$ is $\Sp\colon\SM\to[-1,1]$.
Let $\fN$ be the $V_C$-structure with $\SM^\fN = \{0,1\}$, $\Sp^\fN(0) = -1$ and $\Sp^\fN(1) = 1$, in which every constant of $C$ is interpreted as $0$.
For a variable $\Sx$ of sort $\SM$, the sentence $\supx\Sp(\Sx)$ takes the value $1$ in $\fN$, whereas $(\supx\Sp(\Sx))^{\hat{\fa}} = \sup_n \Sp(\Sc_n^{\SM})^\fa = -1$, the constants $\Sc_n^{\SM}$ being the only closed $V_C$-terms.
Applying the connective $\Sneg$ with $s = -1$ to this sentence reverses the inequality.

Let $\Phi$ be a subset of the set $\Soo^V$ of sentences of $\Loo^V$, and let $\fm\in \prod \Phi$ be given.
We say that the formal type $\fm$ is \textbf{realized} if there is a $V$-structure $\fM$ such that $\Phitp(\fM) = \fm$; in that case, we say that $\fM$ \textbf{realizes} $\fm$.
When $V$ is countable, there are countably many closed $V_C$-terms.
In that case, a $V$-structure can therefore have a canonical expansion to $V_C$, to which Lemma~\ref{lem:henkin-unique} applies, only if its sorts are countable.
The following form of the downward L\"owenheim--Skolem theorem replaces every $V$-structure by one whose sorts are countable and which assigns the same values to the sentences of a given countable fragment.

\begin{lemma}
	\label{lem:downward-ls}
	Assume that the vocabulary $V$ is countable, and let $\frag$ be a countable fragment of $\Loo^V$.
	For every $V$-structure $\fM$ there is a $V$-structure $\fM_0$ whose sorts are countable such that $\phi^{\fM_0} = \phi^\fM$ for every $\phi$ in the set $\SF$ of sentences of $\frag$.
\end{lemma}

\begin{proof}
	The structure $\fM_0$ is obtained by closing countably many elements of $\fM$ under the functions of $\fM$ and under the choice of approximate witnesses for the suprema of formulas of $\frag$.
	Let $\psi\in\frag$, let $\Sy$ be a variable, and let $\bar{x}$ be a finite tuple of distinct variables among which are the free variables of $\supy\psi$.
	For every tuple $\bar{a}\in{}_{\bar{x}}\SM^\fM$ and every integer $k\geq 1$, choose an element $b_{\psi,\Sy,\bar{x},\bar{a},k}\in\sorty^\fM$ such that $\psi^\fM(\bar{a}[b_{\psi,\Sy,\bar{x},\bar{a},k}/\Sy]) > (\supy\psi)^\fM(\bar{a}) - 1/k$, with the notation of the clause for suprema in Section~\ref{sec:lang}.
	For each $\SM\in\MM$, let $A_0^{\SM}$ consist of one element of $\SM^\fM$.
	For each $j<\omega$ and each $\SM\in\MM$, let $A_{j+1}^{\SM}$ consist of the elements of $A_j^{\SM}$, the values in $\SM^\fM$ of the functions $\Sf^\fM$ at tuples from $A_j$, and the witnesses $b_{\psi,\Sy,\bar{x},\bar{a},k}$ in $\SM^\fM$ with $\bar{a}$ from $A_j$.
	Here a tuple is \emph{from $A_j$} if each of its entries lies in the set $A_j^{\SN}$ of its sort $\SN$.
	Since $V$ and $\frag$ are countable, and since there are countably many variables and finite tuples of variables, an induction on $j$ shows that every set $A_j^{\SM}$ is countable; hence so is $A^{\SM} \coloneq \bigcup_j A_j^{\SM}$.
	Moreover, the sets $A_j^{\SM}$ increase with $j$, whence the family $(A^{\SM})_{\SM\in\MM}$ is closed under the functions $\Sf^\fM$.
	Let $\fM_0$ be the $V$-structure whose sorts are the nonempty sets $A^{\SM}$ and whose functions and predicates are the restrictions of those of $\fM$.
	By induction on the rank of $\phi\in\frag$ (see Section~\ref{sec:lang}), we show that $\phi^{\fM_0}(\bar{a}) = \phi^\fM(\bar{a})$ for every tuple $\bar{a}$ from $\fM_0$ of the appropriate sort.
	For atomic $\phi$ this holds because terms take the same values in $\fM_0$ and in $\fM$.
	The steps for $\Sneg$, $\Splus$, $\Sor$ and $\join$ follow from the definitions, the subformulas of members of $\frag$ being members of $\frag$ by (F2).
	Let $\phi$ be of the form $\supy\psi$, and let $\bar{a}$ be a tuple from $\fM_0$ of sort ${}_{\bar{x}}\SM$, where $\bar{x}$ is a finite tuple of distinct variables among which are the free variables of $\phi$.
	The induction hypothesis gives $\phi^{\fM_0}(\bar{a})\leq\phi^\fM(\bar{a})$, since the supremum defining $\phi^{\fM_0}(\bar{a})$ ranges over a subset of $\sorty^\fM$.
	Conversely, the tuple $\bar{a}$ is from $A_j$ for some $j$.
	For every $k\geq 1$, the witness $b\coloneq b_{\psi,\Sy,\bar{x},\bar{a},k}$ is therefore an element of $\fM_0$, and the induction hypothesis applied to $\psi$ gives $\phi^{\fM_0}(\bar{a}) \geq \psi^{\fM_0}(\bar{a}[b/\Sy]) = \psi^{\fM}(\bar{a}[b/\Sy]) > \phi^\fM(\bar{a}) - 1/k$.
	The case of sentences is the lemma.
\end{proof}

Recall from Section~\ref{sec:fragments} that, for a fragment $\frag$ of $\Loo^V$, we write $\frag_C$ for the fragment of $\Loo^C$ generated by $\frag$, and $\SFC$ for its set of sentences.

\begin{proposition}[Henkin Realizability]
	\label{prop:henkin-sat}
	Assume that the vocabulary $V$ is countable, and let $\Phi \subseteq \Soo^V$ and $\fm\in \prod \Phi$ be given.
	\begin{enumerate}[(a)]
		\item The formal type $\fm$ is realized by a $V$-structure whose sorts are countable if and only if it extends to a Henkin type in $\prod \Soo^C$.
		\item If $\Phi\subseteq\SF$ for a countable fragment $\frag$ of $\Loo^V$, then $\fm$ is realized if and only if it extends to an $\SFC$-Henkin type.
	\end{enumerate}
\end{proposition}

\begin{proof}
	The idea is that, since $V$ is countable, a $V$-structure has countable sorts exactly when it has a canonical expansion to $V_C$.
	By Lemma~\ref{lem:henkin-unique}, moreover, the types of canonical $V_C$-structures are exactly the Henkin types.
	In part (b), Lemma~\ref{lem:downward-ls} first replaces the realizing structure by one whose sorts are countable.

	(a) Suppose first that a $V$-structure $\fM$ whose sorts are countable realizes $\fm$.
	Since every sort of $\fM$ is countable and nonempty, there is an expansion $\fM[C]$ of $\fM$ to $V_C$ in which, for each $\SM\in\MM$, the map $n\mapsto(\Sc_n^{\SM})^{\fM[C]}$ carries $\omega$ onto $\SM^\fM$, and every such expansion is canonical.
	By Lemma~\ref{lem:henkin-unique}(a), applied to $\Loo^C$, the full type $\Soo^C\text{-}\!\tp(\fM[C])$ is a Henkin type.
	This Henkin type extends $\fm$, because the value of a $V$-sentence in $\fM[C]$ is its value in $\fM$.
	Conversely, let $\fm'\in\prod\Soo^C$ be a Henkin type that extends $\fm$.
	By Lemma~\ref{lem:henkin-unique}(c), applied to $\Loo^C$, the type $\fm'$ is the full type of a canonical $V_C$-structure $\fN$, whose $V$-reduct therefore realizes $\fm$.
	Since $V_C$ is countable together with $V$, there are only countably many closed $V_C$-terms.
	The sorts of the $V$-reduct of $\fN$ are therefore countable, every element of $\fN$ being the value of such a term.

	(b) If a $V$-structure $\fM$ realizes $\fm$, then so does a $V$-structure $\fM_0$ whose sorts are countable, by Lemma~\ref{lem:downward-ls} and the inclusion $\Phi\subseteq\SF$.
	By Lemma~\ref{lem:henkin-unique}(a), applied to $\frag_C$, the $\SFC$-type of the canonical expansion $\fM_0[C]$ constructed in (a) is an $\SFC$-Henkin type, and it extends $\fm$.
	Conversely, let $\fm'$ be an $\SFC$-Henkin type that extends $\fm$.
	By Lemma~\ref{lem:henkin-unique}(c), applied to $\frag_C$, the type $\fm'$ is the $\SFC$-type of a canonical $V_C$-structure, whose $V$-reduct realizes $\fm$.
\end{proof}

	\section{Topology on the Space of Henkin Types}
	\label{sec:Henkin-topology}

\v{C}ech's form of the Baire category theorem needs a \v{C}ech-complete space, and the Henkin types supply one.
When the vocabulary and the fragment are countable, the Henkin types for the sentences of the expanded fragment form a $G_\delta$ subset of the compact Hausdorff space of all formal types for these sentences, and therefore a \v{C}ech-complete space and a Baire space.

Let $\frag$ be a fragment of $\Loo^V$, and recall from Section~\ref{sec:fragments} that $\SFC$ is the set of sentences of the fragment $\frag_C$ of $\Loo^C$ generated by $\frag$.
We denote the set of $\SFC$-Henkin types by $\HFC$.
As recalled in Section~\ref{sec:topological-prelim}, the product $\prod \SFC = \prod_{\phi\in\SFC}[-r_\phi, r_\phi]$ of Section~\ref{sec:type}, where $r_\phi$ is the a priori bound of $\phi$, is a compact Hausdorff space when endowed with the product topology.
The product topology is generated by the projections (see \cite[Section~2.3 and Proposition~1.4.8]{Engelking:1989}), and the open rays with rational endpoints form a subbasis of the topology of each interval $[-r_\phi, r_\phi]$.
Together, these two facts show that the sets
\begin{align*}
	\open{\synt{\phi > r}} \coloneq \left\{ \fm\in \prod \SFC \colon \phi^\fm > r \right\} \\
	\open{\synt{\phi < r}} \coloneq \left\{ \fm\in \prod \SFC \colon \phi^\fm < r \right\}
\end{align*}
where $\phi\in \SFC$ and $r\in\QQ$, form a subbasis of the product topology.
For an $\SFC$-theory $T$, we write $\open{T}$ for the set of formal types $\fm\in\prod\SFC$ that satisfy every condition of $T$, in the sense of Section~\ref{sec:type}.
We say that $\fm$ \textbf{satisfies} $T$ if $\fm\in\open{T}$; thus each $\open{T}$ is the intersection of the subbasic sets given by the conditions of $T$, and the set $\open{\emptyset}$ is all of $\prod\SFC$.
For each $\phi\in\SFC$, the projection $\fm\mapsto\phi^\fm$ is continuous on $\prod\SFC$ (see \cite[Proposition~1.4.8]{Engelking:1989}).
Hence, for a closed condition $\synt{\phi\geq r}$ with $\phi\in\SFC$, the set $\open{\synt{\phi\geq r}} = \{\fm\in\prod\SFC \colon \phi^\fm\geq r\}$ is closed, being the preimage of the closed ray $[r,\infty)$ under the projection $\fm\mapsto\phi^\fm$.
The sets in the proof of Proposition~\ref{thm:g-delta} below are likewise built from sets defined by equations and inequalities between finitely many coordinates: if $\sigma_1,\dots,\sigma_k\in\SFC$ and $g\colon\RR^k\to\RR$ is continuous, then the set $\{\fm \colon g(\sigma_1^\fm,\dots,\sigma_k^\fm) > 0\}$ is open.
Under the same hypotheses, the set $\{\fm \colon g(\sigma_1^\fm,\dots,\sigma_k^\fm) = 0\}$, being the intersection of the open sets $\{\fm \colon |g(\sigma_1^\fm,\dots,\sigma_k^\fm)| < 1/n\}$ with $n\in\ZZ_+$, is a $G_\delta$ set.

\begin{proposition}
	\label{thm:g-delta}
	If the vocabulary $V$ and the fragment $\frag$ are countable, then the set $\HFC$ of $\SFC$-Henkin types is a $G_\delta$ subset of $\prod \SFC$.
\end{proposition}

\begin{proof}
	The idea is that the Henkin conditions consist of countably many equations, each either between finitely many coordinates or between one coordinate and the supremum of countably many others.
	Each such equation is, moreover, equivalent to countably many requirements defining $G_\delta$ sets.
	All the index sets below are countable: since $V$ is countable, so are $V_C$, the set of variables and the set of closed $V_C$-terms.
	The fragment $\frag_C$ is countable as well, by Section~\ref{sec:fragments}.

	By the continuity of the projections, every set in the following list is a $G_\delta$ subset of $\prod\SFC$.
	In addition, for each $i\in\{1,\dots,6\}$, a formal type belongs to all the sets listed under (H$i$) if and only if it satisfies the Henkin condition (H$i$):
	\begin{enumerate}[(H1)]
		\item the set $\{\fm \colon \unity^\fm = 1\}$;
		\item for $\phi\in\SFC$ and $s\in\QQ$, the set $\{\fm \colon (\Sneg\phi)^\fm = s\,\phi^\fm\}$;
		\item for $\phi,\psi\in\SFC$, the set $\{\fm\colon (\phi\Splus\psi)^\fm=\phi^\fm+\psi^\fm\}$;
		\item for $\phi, \psi\in\SFC$, the set $\{\fm \colon (\phi\Sor\psi)^\fm = \max\{\phi^\fm,\psi^\fm\}\}$;
		\item for each variable $\Sx$ and each formula $\phi\in\frag_C$ whose free variables are among $\Sx$, the closed sets $\{\fm \colon (\supx\phi)^\fm \geq \phi(\tau)^\fm\}$ with $\tau$ a closed $V_C$-term of sort $\sortx$, and the set
		\[
		\bigcap_{n\in\ZZ_+}\bigcup_{\tau}\left\{\fm \colon (\supx\phi)^\fm < \phi(\tau)^\fm + \tfrac{1}{n}\right\},
		\]
		which is a $G_\delta$ set because each union is open;
		\item for each sentence $\join\Phi\in\SFC$, the closed sets $\{\fm \colon (\join\Phi)^\fm \geq \phi^\fm\}$ with $\phi\in\Phi$, and the set
		\[
		\bigcap_{n\in\ZZ_+}\bigcup_{\phi\in\Phi}\left\{\fm \colon (\join\Phi)^\fm < \phi^\fm + \tfrac{1}{n}\right\}.
		\]
	\end{enumerate}
	Each set listed under (H1), (H2), (H3) and (H4), and each closed set listed under (H5) and (H6), is the zero set of a continuous function of finitely many coordinates (such as $\fm\mapsto\max\{\phi^\fm-(\join\Phi)^\fm, 0\}$).
	Each union listed is a union of sets of the form $\{\fm \colon g(\sigma_1^\fm,\dots,\sigma_k^\fm) > 0\}$; in both cases, the remark preceding the proposition therefore applies.
	Being indexed by countable sets, the family of all these sets is countable, and its intersection is $\HFC$.
	Since a countable intersection of $G_\delta$ sets is a $G_\delta$ set, the proposition follows.
\end{proof}

\begin{corollary}
	\label{thm:cech-complete}
	If the vocabulary $V$ and the fragment $\frag$ are countable, then $\HFC$ is \v{C}ech-complete, Polish, and a Baire space.
\end{corollary}

\begin{proof}
	By Proposition~\ref{thm:g-delta}, the space $\HFC$ is a $G_\delta$ subset of the compact Hausdorff space $\prod\SFC$.
	As recalled in Section~\ref{sec:topological-prelim}, every such subset is \v{C}ech-complete, and the Baire category theorem for \v{C}ech-complete spaces (see Theorem~\ref{thm:baire-category}) shows that it is a Baire space.
	Since $\SFC$ is countable, the product $\prod\SFC$ is compact and metrizable, and a $G_\delta$ subset of a Polish space is Polish~\cite[Theorem~3.11]{Kechris:1995}.
\end{proof}

In Section~\ref{sec:omit-type} we express in terms of conditions that a set of Henkin types has empty interior.
For this purpose we use a basis of the topology of $\HFC$ described by finite theories.
Since $\open{S}$ is a finite intersection of subbasic open sets whenever $S$ is a finite $\SFC$-theory, the sets $\HFC\cap\open{S}$ with $S$ finite form a basis of the topology of $\HFC$ (see \cite[Sections~1.1 and~2.1]{Engelking:1989}).

	\section{Omitting Types Theorem}
	\label{sec:omit-type}

The Henkin types satisfying a theory are the types of the canonical structures (those in which every element is the value of a closed term) whose reducts are models of the theory.
To omit countably many types in a model of the theory, it therefore suffices to find such a Henkin type that realizes none of them.
In this section we show that, when the vocabulary and the fragment are countable and the theory has a model, the Henkin types satisfying the theory form a nonempty Baire space.
We also show that the required Henkin type exists if, for each instance of each of the given types, the members of this space satisfying the instance form a nowhere dense set.
The members of this space that realize one of the given types then form a meager subset, and the Baire category theorem provides a Henkin type outside it.

Henceforth, assume that the vocabulary $V$ and the fragment $\frag$ of Section~\ref{sec:Henkin-topology} are countable.
We also fix an $\SF$-theory $B$ (a theory whose conditions are on sentences of $\frag$) that has a model.
Recall from Section~\ref{sec:Henkin-topology} that $\HFC$ is the set of $\SFC$-Henkin types, where $\SFC$ is the set of sentences of the fragment $\frag_C$ of $\Loo^C$ generated by $\frag$.
The set $\HFC\cap\open{B}$ of $\SFC$-Henkin types that satisfy $B$ (see Section~\ref{sec:Henkin-topology}) is nonempty.
Indeed, by Lemma~\ref{lem:downward-ls}, the theory $B$ has a model $\fM$ whose sorts are countable.
The canonical expansion $\fM[C]$ of $\fM$ constructed in the proof of Proposition~\ref{prop:henkin-sat}(a) satisfies $B$, since the values of $V$-sentences in $\fM[C]$ are their values in $\fM$.
By Lemma~\ref{lem:henkin-unique}(a), the $\SFC$-type of $\fM[C]$ belongs to $\HFC\cap\open{B}$.
Conversely, by Lemma~\ref{lem:henkin-unique}(c), every member of $\HFC\cap\open{B}$ is the $\SFC$-type of a canonical $V_C$-structure, whose $V$-reduct is a model of $B$.
Thus the members of $\HFC\cap\open{B}$ are exactly the $\SFC$-types of the canonical $V_C$-structures whose $V$-reducts are models of $B$.

Since $\SF$ is countable, so is the theory $B$, whence the set $\open{B}$ is a countable intersection of open sets.
By Proposition~\ref{thm:g-delta}, the set $\HFC\cap\open{B}$ is therefore a $G_\delta$ subset of the compact Hausdorff space $\prod\SFC$.
As recalled in Section~\ref{sec:topological-prelim}, a $G_\delta$ subset of a compact Hausdorff space is \v{C}ech-complete; by Theorem~\ref{thm:baire-category}, the space $\HFC\cap\open{B}$ is a nonempty Baire space.

Let $\bar{x}$ be a finite tuple of distinct variables.
An \textbf{$\frag$-type in $\bar{x}$} is a set $\Gamma(\bar{x})$ of expressions $\synt{\phi(\bar{x})\geq r}$, where $r\in\QQ$ and $\phi\in\frag$ is a formula whose free variables are among $\bar{x}$.
Substituting closed terms for $\bar{x}$ turns each of these expressions into a closed condition in the sense of Section~\ref{sec:cond-theories}.
For a tuple $\bar{\tau}$ of closed $V_C$-terms of sort ${}_{\bar{x}}\SM$, the \textbf{instance} $\Gamma(\bar{\tau})$ of $\Gamma(\bar{x})$ is the union of the closed conditions $\synt{\phi(\bar{\tau})\geq r}$ for the members $\synt{\phi(\bar{x})\geq r}$ of $\Gamma(\bar{x})$.
Each of these conditions is a countable theory; the instance $\Gamma(\bar{\tau})$ is an $\SFC$-theory because the fragment $\frag_C$ is closed under substitution of terms (see (F4) in Section~\ref{sec:fragments}).
An $\SFC$-Henkin type $\fm\in \HFC$ \textbf{realizes} $\Gamma(\bar{x})$ if it satisfies some instance $\Gamma(\bar{\tau})$, that is, if $\fm\in\open{\Gamma(\bar{\tau})}$ for some tuple $\bar{\tau}$; otherwise it \textbf{omits} $\Gamma(\bar{x})$.
Similarly, if $\fM$ is a $V$-structure, a tuple $\bar{a}\in{}_{\bar{x}}\SM^\fM$ \textbf{satisfies} $\Gamma(\bar{x})$ in $\fM$ if $\phi^\fM(\bar{a})\geq r$ for every member $\synt{\phi(\bar{x})\geq r}$ of $\Gamma(\bar{x})$.
We say that $\Gamma(\bar{x})$ \textbf{has nonempty interior in $\HFC\cap\open{B}$} if there exist a tuple $\bar{\tau}$ of closed $V_C$-terms of sort ${}_{\bar{x}}\SM$ and a finite $\SFC$-theory $S$ such that
\[
\emptyset \neq \HFC\cap \open{B} \cap \open{S} \subseteq \open{\Gamma(\bar{\tau})}.
\]
Otherwise, we say that $\Gamma(\bar{x})$ \textbf{has empty interior in $\HFC\cap\open{B}$}.
The sets $\HFC\cap\open{B}\cap\open{S}$, with $S$ a finite $\SFC$-theory, form a basis of the topology of $\HFC\cap\open{B}$ (see Section~\ref{sec:Henkin-topology}).
Consequently, the type $\Gamma(\bar{x})$ has nonempty interior in this space exactly when, for some tuple $\bar{\tau}$, the closed subset $\HFC\cap\open{B}\cap\open{\Gamma(\bar{\tau})}$ of this space has nonempty interior.
An $\frag$-type is expressed by closed conditions so that this subset, being an intersection of the closed sets $\open{\synt{\phi(\bar{\tau})\geq r}}$, is closed (see Section~\ref{sec:Henkin-topology}).
When this subset has empty interior, its complement in $\HFC\cap\open{B}$ is therefore a dense open set.
The Baire category theorem applies to countable intersections of such dense open sets, as the proof of Theorem~\ref{thm:omitting-types} shows.

\begin{theorem}[Omitting Types Theorem]
	\label{thm:omitting-types}
	Assume that the vocabulary $V$ and the fragment $\frag$ are countable, and let $B$ be an $\SF$-theory that has a model.
	Let $(\Gamma_n(\bar{x}_n))_{n<\omega}$ be a sequence of $\frag$-types each of which has empty interior in $\HFC\cap\open{B}$.
	Then there is an $\SFC$-Henkin type $\fm\in \HFC\cap\open{B}$ that omits every $\Gamma_n(\bar{x}_n)$, and the theory $B$ has a model whose sorts are countable and in which, for every $n<\omega$, no tuple of elements satisfies $\Gamma_n(\bar{x}_n)$.
\end{theorem}

\begin{proof}
	The idea is to write the set of Henkin types that satisfy $B$ and omit every $\Gamma_n(\bar{x}_n)$ as a countable intersection of open dense subsets of the nonempty Baire space $\HFC\cap\open{B}$.
	Fix $n<\omega$, let $\bar{x}\coloneq\bar{x}_n$ and $\Gamma\coloneq \Gamma_n$, and let $\bar{\tau}$ be a tuple of closed $V_C$-terms of sort ${}_{\bar{x}}\SM$.
	Since the set $\open{\synt{\phi(\bar{\tau})\geq r}}$ is closed for every member $\synt{\phi(\bar{x})\geq r}$ of $\Gamma(\bar{x})$ (see Section~\ref{sec:Henkin-topology}), the set $\HFC\cap\open{B}\cap\open{\Gamma(\bar{\tau})}$ is closed in the space $\HFC\cap\open{B}$.
	By the definition of empty interior, this closed set has empty interior in $\HFC\cap\open{B}$, and the complement
	\[
	O_{\Gamma,\bar{\tau}} \coloneq \{ \fm\in\HFC\cap\open{B} \colon \fm \text{ does not satisfy } \Gamma(\bar{\tau}) \}
	\]
	is therefore an open dense subset of $\HFC\cap\open{B}$.
	Since $V$ is countable, there are only countably many tuples $\bar{\tau}$; consequently, the set $O_\Gamma \coloneq \bigcap_{\bar{\tau}} O_{\Gamma,\bar{\tau}}$ of members of $\HFC\cap\open{B}$ that omit $\Gamma(\bar{x})$ is a countable intersection of open dense subsets of $\HFC\cap\open{B}$.
	As $n$ was arbitrary, the same holds for each set $O_{\Gamma_n}$, and hence for the set $O\coloneq \bigcap_{n<\omega} O_{\Gamma_n}$ of members of $\HFC\cap\open{B}$ that omit $\Gamma_n(\bar{x}_n)$ for every $n<\omega$.
	Since $\HFC\cap \open{B}$ is a nonempty Baire space, the set $O$ is dense in it, and in particular nonempty.

	Let $\fm$ be a member of $O$.
	By Lemma~\ref{lem:henkin-unique}(c), the type $\fm$ is the $\SFC$-type of a canonical $V_C$-structure $\fN$, whose $V$-reduct $\fM$ is a model of $B$.
	The sorts of $\fM$ are countable, since every element of $\fN$ is the value of one of the countably many closed $V_C$-terms.
	Now fix $n<\omega$ and a tuple $\bar{a}\in{}_{\bar{x}_n}\SM^\fM$.
	Since $\fN$ is canonical, we may write $\bar{a} = \bar{\tau}^\fN$ for a tuple $\bar{\tau}$ of closed $V_C$-terms.
	By the substitution lemma (Lemma~\ref{lem:substitution}), with $\bar{z}$ empty, we have $\phi^\fM(\bar{a}) = \phi^\fN(\bar{\tau}^\fN) = \phi(\bar{\tau})^\fm$ for every formula $\phi\in\frag$ whose free variables are among $\bar{x}_n$.
	Since $\fm$ omits $\Gamma_n(\bar{x}_n)$, it does not satisfy $\Gamma_n(\bar{\tau})$, and the tuple $\bar{a}$ therefore does not satisfy $\Gamma_n(\bar{x}_n)$ in $\fM$.
\end{proof}

\part{Scott's Theorem and Its Topological Form}

\section{Scott Sentences}
\label{sec:scott}

In 1965, Scott proved that every countable structure for a countable vocabulary is characterized, among countable structures and up to isomorphism, by a single sentence of the infinitary logic of two-valued structures~\cite{Scott:1965}.
Ben Yaacov, Doucha, Nies and Tsankov extended the theorem to separable metric structures whose predicates are uniformly continuous with respect to a complete metric, with moduli prescribed by the signature~\cite{BenYaacov-Doucha-Nies-Tsankov:2017}.
In this section we extend Scott's theorem to the general structures of this paper, which require no prescribed metric and no prescribed continuity.
In $\Loo$, the value of a sentence is invariant under approximate back-and-forth equivalence, so that a sentence registers a structure only up to approximation.
Accordingly, the theorem characterizes every structure whose sorts are countable, among all structures, up to approximate back-and-forth equivalence, and, for structures with a predicate for equality in which every atomic formula $\psi$ takes only the two values $r_\psi$ and $-r_\psi$, the theorem yields the content of Scott's theorem.

Throughout the section, the vocabulary $V=(\MM,F,P)$ is countable.
We use the following abbreviations for formulas $\phi$ and $\psi$ of $\Loo^V$, a variable $\Sx$ and a sequence $\Phi=(\phi_n)_{n<\omega}$ whose join is a formula: $\synt{-\phi}$ is $\Sneg\phi$ with $s=-1$; $\phi\Sand\psi$ is $\synt{-}(\synt{-\phi}\Sor\synt{-\psi})$; $\Sinf{\Sx}\phi$ is $\synt{-}\supx\synt{-\phi}$; $\meet\Phi$ is $\synt{-}\join_n\synt{-\phi_n}$; $\synt{|\phi|}$ is $\phi\Sor\synt{-\phi}$; and $\phi\Sdot\psi$ is $(\phi\Splus\synt{-\psi})\Sor(\synt{0\cdot}\unity)$.
Apart from $\synt{-\phi}$, their values are $\min\{\phi^\fM,\psi^\fM\}$, the infimum of the values of $\phi$ over the sort of $\Sx$, the infimum of the values of the members of $\Phi$, $|\phi^\fM|$ and $\max\{\phi^\fM-\psi^\fM,0\}$, and their a priori bounds are $\max\{r_\phi,r_\psi\}$, $r_\phi$, $\sup_nr_{\phi_n}$, $r_\phi$ and $r_\phi+r_\psi$, respectively.

Real numbers enter the language as countable meets of rational multiples of the logical unity constant.

\begin{lemma}
\label{lem:constants}
Let real numbers $t$ and $\epsilon>0$ be given.
Then there is a $V$-sentence $\underline{t}$ of $\Loo^V$ that takes the value $t$ in every $V$-structure and has a priori bound at most $|t|+\epsilon$.
\end{lemma}

\begin{proof}
Choose rational numbers $q_n$ with $t\leq q_n\leq t+\epsilon$ and $q_n\to t$, and let $\underline{t}$ be the meet of the sentences $\synt{q_n\cdot}\unity$.
These sentences have the common a priori bound $|t|+\epsilon$, and the value of their meet is $\inf_nq_n=t$.
\end{proof}

For each finite sequence $(\SM_i)_{i<k}$ of sort descriptors, we use the tuple $\bar{x}=(\Sx_i^{\SM_i})_{i<k}$ of variables, whose entries are distinct because their indices are; a tuple of elements of a structure of sort ${}_{\bar{x}}\SM$ is said to be of \textbf{sort} $(\SM_i)_{i<k}$.
The variable that extends such a tuple by an entry of sort $\SM$ is $\Sx_k^{\SM}$, distinct from the entries of $\bar{x}$.
An \textbf{atomic formula in $\bar{x}$} is an atomic $V$-formula whose variables are among $\bar{x}$; since $V$ is countable, there are countably many of them.

\begin{definition}
\label{def:system}
Let $V$-structures $\fM$ and $\fN$ and a real number $r>0$ be given.
An \textbf{$r$-back-and-forth system} between $\fM$ and $\fN$ is a set $I$ of pairs $(\bar{a},\bar{b})$ of tuples of the same sort, $\bar{a}$ of elements of $\fM$ and $\bar{b}$ of elements of $\fN$, such that the pair of empty tuples belongs to $I$ and, for every pair $(\bar{a},\bar{b})$ in $I$ of sort $(\SM_i)_{i<k}$, with $\bar{x}=(\Sx_i^{\SM_i})_{i<k}$:
\begin{enumerate}[(i)]
\item $|\psi^\fM(\bar{a})-\psi^\fN(\bar{b})|\leq r\,r_\psi$ for every atomic formula $\psi$ in $\bar{x}$;
\item for every sort descriptor $\SM\in\MM$ and every $c\in\SM^\fM$ there is $d\in\SM^\fN$ with $(\bar{a}c,\bar{b}d)\in I$;
\item for every sort descriptor $\SM\in\MM$ and every $d\in\SM^\fN$ there is $c\in\SM^\fM$ with $(\bar{a}c,\bar{b}d)\in I$.
\end{enumerate}
The structures $\fM$ and $\fN$ are \textbf{approximately equivalent} if there is an $r$-back-and-forth system between them for every real number $r>0$.
\end{definition}

\begin{lemma}
\label{lem:transfer}
Let $V$-structures $\fM$ and $\fN$, a real number $r>0$ and an $r$-back-and-forth system $I$ between them be given.
Then $|\phi^\fM(\bar{a})-\phi^\fN(\bar{b})|\leq r\,r_\phi$ for every pair $(\bar{a},\bar{b})$ in $I$, of sort $(\SM_i)_{i<k}$ say, and every formula $\phi(\bar{x})$ of $\Loo^V$, where $\bar{x}=(\Sx_i^{\SM_i})_{i<k}$.
\end{lemma}

\begin{proof}
We argue by induction on the rank of $\phi$ (see Section~\ref{sec:lang}), for all pairs in $I$ at once.
For atomic formulas the inequality is condition~(i).
For $\Sneg\phi$, $\phi\Splus\psi$ and $\phi\Sor\psi$, the inequality follows from the inequalities for $\phi$ and $\psi$, since $|\max\{u,v\}-\max\{u',v'\}|\leq\max\{|u-u'|,|v-v'|\}$ for real numbers $u$, $v$, $u'$ and $v'$; for $\join\Phi$, the inequality follows from $|\sup_nu_n-\sup_nu'_n|\leq\sup_n|u_n-u'_n|$ for bounded sequences of real numbers.
Let $\phi(\bar{x})$ be of the form $\supy\psi$, with $\Sy$ of sort $\SM$.
Rename the bound variables of $\psi$ away from the variable $\Sx_k^{\SM}$, and let $\psi'(\bar{x}\Sx_k^{\SM})$ be obtained by substituting $\Sx_k^{\SM}$ for the free occurrences of $\Sy$.
The substitution preserves the rank, the a priori bound and the admissibility of every join, and the formulas $\supy\psi$ and $\synt{\sup_{\Sx_k^{\SM}}}\psi'$ take the same values; the induction hypothesis therefore applies to $\psi'$.
Let $c\in\SM^\fM$ be given; condition~(ii) gives $d\in\SM^\fN$ with $(\bar{a}c,\bar{b}d)\in I$, and the induction hypothesis gives
\[
\psi'^{\fM}(\bar{a}c)\leq\psi'^{\fN}(\bar{b}d)+r\,r_\psi\leq\phi^\fN(\bar{b})+r\,r_\phi.
\]
Taking the supremum over $c$ gives $\phi^\fM(\bar{a})\leq\phi^\fN(\bar{b})+r\,r_\phi$, and condition~(iii) gives the reverse inequality in the same way.
\end{proof}

\begin{corollary}
\label{cor:equivalence}
Approximately equivalent $V$-structures have the same full type.
\end{corollary}

\begin{proof}
Applied to the pair of empty tuples, Lemma~\ref{lem:transfer} gives $|\phi^\fM-\phi^\fN|\leq r\,r_\phi$ for every sentence $\phi$ of $\Loo^V$ and every $r>0$.
\end{proof}

Let a $V$-structure $\fM$ whose sorts are countable be given.
The Scott formulas of $\fM$, defined below, measure the error with which a tuple of another structure matches a tuple of $\fM$ in the back-and-forth game.
For an atomic formula $\psi$ in $\bar{x}=(\Sx_i^{\SM_i})_{i<k}$ with $r_\psi>0$, choose a rational number $c_\psi$ with $1/(6r_\psi)\leq c_\psi\leq1/(3r_\psi)$ and a rational number $\epsilon_\psi$ with $0<\epsilon_\psi\leq r_\psi$.
For a tuple $\bar{a}$ of elements of $\fM$ of sort $(\SM_i)_{i<k}$, let $t=\psi^\fM(\bar{a})$, and let
\[
\delta_{\psi,\bar{a}}(\bar{x})\coloneq\synt{c_\psi\cdot}\,\synt{|}\psi\Splus\synt{-}\underline{t}\,\synt{|},
\]
where $\underline{t}$ is the sentence of Lemma~\ref{lem:constants} for $t$ and $\epsilon_\psi$.
The a priori bound of $\delta_{\psi,\bar{a}}$ is at most $c_\psi(r_\psi+|t|+\epsilon_\psi)\leq3c_\psi r_\psi\leq1$, and, for every $V$-structure $\fN$ and every tuple $\bar{b}$ of elements of $\fN$ of sort $(\SM_i)_{i<k}$,
\begin{equation}
\label{eq:atomic}
|\psi^\fN(\bar{b})-\psi^\fM(\bar{a})|\leq6r_\psi\,\delta_{\psi,\bar{a}}^\fN(\bar{b}).
\end{equation}

\begin{definition}
\label{def:scott-formulas}
For every countable ordinal $\alpha$ and every tuple $\bar{a}$ of elements of $\fM$ of sort $(\SM_i)_{i<k}$, the \textbf{Scott formula} $\phi^\alpha_{\bar{a}}(\bar{x})$, with $\bar{x}=(\Sx_i^{\SM_i})_{i<k}$, is defined by recursion on $\alpha$.
The formula $\phi^0_{\bar{a}}$ is the join of the formula $\synt{0\cdot}\unity$ and of the formulas $\delta_{\psi,\bar{a}}$ for the atomic formulas $\psi$ in $\bar{x}$ with $r_\psi>0$.
For every countable ordinal $\alpha$, the formula $\phi^{\alpha+1}_{\bar{a}}$ is the join of the formula $\phi^\alpha_{\bar{a}}$ and, for every sort descriptor $\SM\in\MM$, of the formulas
\[
\Sinf{\Sx_k^{\SM}}\phi^\alpha_{\bar{a}c}(\bar{x}\Sx_k^{\SM})\quad\text{for }c\in\SM^\fM,
\qquad
\synt{\sup_{\Sx_k^{\SM}}}\meet_{c\in\SM^\fM}\phi^\alpha_{\bar{a}c}(\bar{x}\Sx_k^{\SM}).
\]
For every countable limit ordinal $\lambda$, the formula $\phi^\lambda_{\bar{a}}$ is the join of the formulas $\phi^\alpha_{\bar{a}}$ with $\alpha<\lambda$.
\end{definition}

By induction on $\alpha$, every Scott formula has a priori bound at most $1$, so that every join and every meet in the recursion is a formula of $\Loo^V$, the sorts of $\fM$, the set $\MM$, the set of atomic formulas in $\bar{x}$ and the set of ordinals below a countable ordinal being countable.
By the same induction, every Scott formula takes values in $[0,1]$, and $\phi^\alpha_{\bar{a}}$ is at most $\phi^\beta_{\bar{a}}$ in every $V$-structure for $\alpha\leq\beta$.
In the recursion, the infima measure the moves from $\fM$ to the other structure and the suprema the moves back.

\begin{lemma}
\label{lem:self-match}
The value of every Scott formula $\phi^\alpha_{\bar{a}}$ at the tuple $\bar{a}$ in $\fM$ is $0$.
\end{lemma}

\begin{proof}
For $\alpha=0$, every formula $\delta_{\psi,\bar{a}}$ vanishes at $\bar{a}$.
For a successor ordinal $\alpha+1$, let $\SM\in\MM$ be given.
The induction hypothesis gives $\inf_{b}\phi^{\alpha,\fM}_{\bar{a}c}(\bar{a}b)\leq\phi^{\alpha,\fM}_{\bar{a}c}(\bar{a}c)=0$ for every $c\in\SM^\fM$, and, for every $e\in\SM^\fM$, the infimum of the values $\phi^{\alpha,\fM}_{\bar{a}c}(\bar{a}e)$ over $c\in\SM^\fM$ is at most $\phi^{\alpha,\fM}_{\bar{a}e}(\bar{a}e)=0$.
The limit case is immediate.
\end{proof}

\begin{lemma}
\label{lem:stable}
There is a countable ordinal $\alpha_\fM$ such that $\phi^{\alpha_\fM+1,\fM}_{\bar{a}}(\bar{b})=\phi^{\alpha_\fM,\fM}_{\bar{a}}(\bar{b})$ for all tuples $\bar{a}$ and $\bar{b}$ of elements of $\fM$ of the same sort.
\end{lemma}

\begin{proof}
For fixed tuples $\bar{a}$ and $\bar{b}$, the values $\phi^{\alpha,\fM}_{\bar{a}}(\bar{b})$, for $\alpha<\omega_1$, form a nondecreasing family of real numbers.
The open intervals between the values at $\alpha$ and at $\alpha+1$, for the ordinals $\alpha$ at which the family increases, are pairwise disjoint, so these ordinals are countably many; since $\omega_1$ is regular, there is a countable ordinal from which the value at every ordinal equals the value at its successor.
At a countable limit ordinal $\lambda$ the value is the supremum of the earlier values, since $\phi^\lambda_{\bar{a}}$ is the join of the earlier formulas;
by induction, the family is constant from that ordinal on.
There are countably many pairs $(\bar{a},\bar{b})$, and the supremum of countably many countable ordinals is countable.
\end{proof}

\begin{definition}
\label{def:scott-sentence}
The \textbf{Scott sentence} $\sigma_\fM$ of $\fM$ is the join of the sentence $\phi^{\alpha_\fM+1}_{\varnothing}$, for the empty tuple $\varnothing$, and of the sentences
\[
\synt{\sup_{\Sx_0^{\SM_0}}}\cdots\synt{\sup_{\Sx_{k-1}^{\SM_{k-1}}}}\bigl(\phi^{\alpha_\fM+1}_{\bar{a}}\Sdot\phi^{\alpha_\fM}_{\bar{a}}\bigr)
\]
for the tuples $\bar{a}$ of elements of $\fM$, of every sort $(\SM_i)_{i<k}$, where the ordinal $\alpha_\fM$ is given by Lemma~\ref{lem:stable}.
\end{definition}

The members of this join have a priori bound at most $2$, so $\sigma_\fM$ is a sentence of $\Loo^V$; the values of $\sigma_\fM$ are nonnegative, and the value of $\sigma_\fM$ in $\fM$ is $0$, by Lemmas~\ref{lem:self-match} and~\ref{lem:stable}.

\begin{theorem}[Scott's theorem for infinitary real-valued logic]
\label{thm:scott}
Assume that the vocabulary $V$ is countable, let $\fM$ be a $V$-structure whose sorts are countable, and let $\sigma_\fM$ be the Scott sentence of $\fM$.
For every $V$-structure $\fN$, the following conditions are equivalent:
\begin{enumerate}[(a)]
\item $\fN$ satisfies the closed condition $\synt{\sigma_\fM\leq 0}$;
\item $\fM$ and $\fN$ are approximately equivalent;
\item $\fM$ and $\fN$ have the same full type.
\end{enumerate}
\end{theorem}

\begin{proof}
Condition (b) implies (c) by Corollary~\ref{cor:equivalence}, and (c) implies (a) because $\sigma_\fM^\fM=0$.

Suppose (a); since $\sigma_\fM$ takes nonnegative values, $\sigma_\fM^\fN=0$.
Let a real number $r>0$ be given, and let $I$ be the set of pairs $(\bar{a},\bar{b})$ of tuples of the same sort, $\bar{a}$ of elements of $\fM$ and $\bar{b}$ of elements of $\fN$, such that $\phi^{\alpha_\fM,\fN}_{\bar{a}}(\bar{b})<r/6$.
Since $\sigma_\fM^\fN=0$, the value in $\fN$ of $\phi^{\alpha_\fM+1}_{\varnothing}$ is $0$, and so is the smaller value of $\phi^{\alpha_\fM}_{\varnothing}$; the pair of empty tuples belongs to $I$.
Moreover, $\phi^{\alpha_\fM+1,\fN}_{\bar{a}}(\bar{b})\leq\phi^{\alpha_\fM,\fN}_{\bar{a}}(\bar{b})$ for all tuples $\bar{a}$ of $\fM$ and $\bar{b}$ of $\fN$ of the same sort.
Let $(\bar{a},\bar{b})$ be a pair in $I$, of sort $(\SM_i)_{i<k}$.
For condition~(i) of Definition~\ref{def:system}, let $\psi$ be an atomic formula in $\bar{x}=(\Sx_i^{\SM_i})_{i<k}$; if $r_\psi=0$, both values of $\psi$ vanish, and otherwise $\phi^{\alpha_\fM}_{\bar{a}}\geq\phi^0_{\bar{a}}\geq\delta_{\psi,\bar{a}}$ and \eqref{eq:atomic} give $|\psi^\fN(\bar{b})-\psi^\fM(\bar{a})|<6r_\psi\cdot r/6=r\,r_\psi$.
For condition~(ii), let $\SM\in\MM$ and $c\in\SM^\fM$ be given.
Then $\inf_{d\in\SM^\fN}\phi^{\alpha_\fM,\fN}_{\bar{a}c}(\bar{b}d)\leq\phi^{\alpha_\fM+1,\fN}_{\bar{a}}(\bar{b})\leq\phi^{\alpha_\fM,\fN}_{\bar{a}}(\bar{b})<r/6$, so there is $d\in\SM^\fN$ with $(\bar{a}c,\bar{b}d)\in I$.
For condition~(iii), let $\SM\in\MM$ and $d\in\SM^\fN$ be given.
Then the infimum of the values $\phi^{\alpha_\fM,\fN}_{\bar{a}c}(\bar{b}d)$ over $c\in\SM^\fM$ is at most $\phi^{\alpha_\fM+1,\fN}_{\bar{a}}(\bar{b})<r/6$, so there is $c\in\SM^\fM$ with $(\bar{a}c,\bar{b}d)\in I$.
Thus $I$ is an $r$-back-and-forth system, and (b) holds.
\end{proof}

The proof rests on a strict inequality: a match with error less than $r/6$ extends, one element at a time, to matches with error less than $r/6$, so that the error does not accumulate along the back-and-forth game.
The structure $\fN$ is arbitrary, and Theorem~\ref{thm:scott} is the real-valued counterpart of Scott's theorem together with the characterization of infinitary equivalence by back-and-forth systems.
For structures with a predicate for equality in which every atomic formula $\psi$ takes only the two values $r_\psi$ and $-r_\psi$, approximate equivalence among countable structures is isomorphism (Corollary~\ref{cor:scott-classical}).

\begin{corollary}[The content of Scott's theorem]
\label{cor:scott-classical}
Assume that $V$ has a single sort descriptor $\SM$ and a predicate symbol $\symb{E}\colon\SM\times\SM\to[-1,1]$.
Let $\fM$ and $\fN$ be $V$-structures whose sorts are countable, in which every atomic formula $\psi$ other than $\unity$ takes only the two values $r_\psi$ and $-r_\psi$ and $\symb{E}$ takes the value $1$ exactly on the diagonal.
Then $\fN$ satisfies $\synt{\sigma_\fM\leq0}$ if and only if $\fM$ and $\fN$ are isomorphic.
\end{corollary}

\begin{proof}
Isomorphic structures have the same full type, and $\sigma_\fM^\fM=0$; so isomorphism gives $\sigma_\fM^\fN=0$.
Conversely, suppose that $\sigma_\fM^\fN=0$; by Theorem~\ref{thm:scott} there is an $r$-back-and-forth system $I$ between $\fM$ and $\fN$ for $r=1$.
The two values $r_\psi$ and $-r_\psi$ of an atomic formula $\psi$ differ by $2r_\psi>r\,r_\psi$, so condition~(i) of Definition~\ref{def:system} gives $\psi^\fM(\bar{a})=\psi^\fN(\bar{b})$ for every pair $(\bar{a},\bar{b})$ in $I$ and every atomic formula $\psi$ in the corresponding variables.
Applied to the atomic formulas $\symb{E}(\Sx_i,\Sx_j)$ and $\symb{E}(\Sf(\bar{\tau}),\Sx_j)$, for function symbols $\Sf$ and tuples $\bar{\tau}$ of variables, this shows that the map $a_i\mapsto b_i$ is well defined, injective, and compatible with the predicates and with the functions on the entries of the tuples.
Enumerate the sorts of $\fM$ and $\fN$, and choose pairs in $I$ alternately by conditions~(ii) and~(iii), each extending the previous one, so that every element of $\fM$ appears in the first tuples and every element of $\fN$ in the second; the union of the corresponding maps is an isomorphism.
\end{proof}

Approximation is a feature of the logic itself, as the following example shows.

\begin{example}
\label{ex:dense}
Assume that $V$ has a single sort descriptor $\SM$, no function symbols, and, besides $\unity$, a single predicate symbol $\Sp\colon\SM\to[-1,1]$.
Let $\fM$ and $\fN$ be the $V$-structures with $\SM^\fM=\{0\}\cup\{1/n\colon n\geq1\}$ and $\SM^\fN=\{1/n\colon n\geq1\}$, in which $\Sp$ is interpreted as the identity.
For $r>0$, let $I$ consist of the pairs $(\bar{a},\bar{b})$ of tuples of the same length with $|a_i-b_i|\leq r$ for every $i$.
The atomic formulas in $\bar{x}$ other than $\unity$ are the formulas $\Sp(\Sx_i)$, so condition~(i) of Definition~\ref{def:system} holds; condition~(ii) holds because every element of $\fM$ lies within $r$ of an element of $\fN$, and condition~(iii) because every element of $\fN$ is an element of $\fM$.
Thus $\fM$ and $\fN$ are approximately equivalent, and, by Corollary~\ref{cor:equivalence}, they have the same full type.
They are not isomorphic, since $\Sp$ vanishes at an element of $\fM$ and at no element of $\fN$.
The formula $\join_n\bigl(\synt{n\cdot}\Sp(\Sx)\Sand\unity\bigr)$, which would separate the element $0$ of $\fM$ from the elements of $\fN$, is not a formula of $\Loo^V$, since the a priori bounds $n$ of the members of the join have no common bound.
\end{example}

\begin{remark}
\label{rem:completions}
The formulas of $\Loo^V$ induce on each $V$-structure $\fN$ a pseudometric on each sort: for elements $a$ and $a'$ of $\SM^\fN$, the number $d_\fN(a,a')$ is the supremum of $|\phi^\fN(\bar{b}a)-\phi^\fN(\bar{b}a')|/r_\phi$ over the formulas $\phi$ of $\Loo^V$ with $r_\phi>0$ and the tuples $\bar{b}$ of elements of $\fN$ of the appropriate sort.
By definition, every predicate of $\fN$ is Lipschitz for these pseudometrics, with constant given by the a priori bound of the predicate symbol in each argument.
In Example~\ref{ex:dense}, $d_\fN(a,a')=|a-a'|$ for distinct elements $a$ and $a'$: the formula $\Sp(\Sx)$ gives one inequality, and Lemma~\ref{lem:transfer}, applied to the $s$-back-and-forth system between $\fN$ and itself consisting of the pairs $(\bar{c},\bar{e})$ with $|c_i-e_i|\leq s$ for every $i$, where $s=|a-a'|$, gives the other.
The completion of $\fN$ for $d_\fN$, with $\Sp$ extended by continuity, is therefore a copy of $\fM$.
Keisler obtained a distance predicate from first-order formulas for every countable vocabulary~\cite[Section~3.3]{Keisler:2023}.
Approximately equivalent structures need not, however, have isomorphic completions for these pseudometrics, even when both structures are already complete, as the next example shows.
\end{remark}

\begin{example}
\label{ex:complete}
Assume that $V$ has a single sort descriptor $\SM$, no function symbols, and, besides $\unity$, predicate symbols $\Sp\colon\SM\to[-1,1]$ and $\symb{E}\colon\SM\times\SM\to[-1,1]$.
Let $\fM$ and $\fN$ be the $V$-structures with $\SM^\fM=\QQ\cap[0,1]$ and $\SM^\fN=\QQ\cap(0,1]$, in which $\Sp$ is interpreted as the identity and $\symb{E}$ takes the value $1$ on the diagonal and $-1$ elsewhere.
For $r>0$, the pairs $(\bar{a},\bar{b})$ of tuples of the same length with $a_i=a_j$ exactly when $b_i=b_j$ and with $|a_i-b_i|<r$ for every $i$ form an $r$-back-and-forth system, since every interval of length $2r$ centered at a point of $[0,1]$ contains infinitely many elements of each structure; so $\fM$ and $\fN$ are approximately equivalent.
For each of the two structures and distinct elements $a$ and $a'$ of it, the formula $\symb{E}(\Sy,\Sx)$ with the parameter $a$ for $\Sy$ gives $d(a,a')\geq2$, where $d$ is the pseudometric of that structure (Remark~\ref{rem:completions}), and the a priori bounds give $d(a,a')\leq2$; so both pseudometrics are discrete and both structures are complete.
The structures are not isomorphic, since $\Sp$ takes the value $0$ in $\fM$ and not in $\fN$.
\end{example}

\section{Scott's Theorem in the Space of Countable Structures}
\label{sec:topological-scott}

In classical logic, Scott's theorem has a topological form.
The structures whose universe is the set of natural numbers form a Polish space, the isomorphism classes are the orbits of the group of permutations of the natural numbers, and Scott's theorem implies that every orbit is a Borel set.
López-Escobar's theorem identifies the Borel sets that are invariant under the group with the sets defined by sentences of the infinitary logic~\cite{LopezEscobar:1965,Vaught:1974}.

We give Theorem~\ref{thm:scott} a corresponding topological form.
Since general structures carry no equality, we replace the group of permutations by the monoid of surjections of the natural numbers, acting on atomic diagrams.
For structures on the natural numbers, exact equivalence means that the orbits of their diagrams meet, and approximate equivalence means that these orbits come arbitrarily close for a uniform distance between diagrams.
The resulting theorem (Theorem~\ref{thm:topological-scott}) involves only a Polish space, the action of a Polish monoid and a distance, and the logic enters only the proof, through the Scott sentences of Section~\ref{sec:scott}.

A version of López-Escobar's theorem then marks the reach of this topological form: the bounded Borel functions invariant under exact equivalence are the values of the sentences of a richer logic, and some of them separate approximately equivalent structures, so that replacing the functions of Theorem~\ref{thm:topological-scott} by all such functions changes the equivalence they detect, from approximate equivalence to exact equivalence.

Throughout the section, the vocabulary $V=(\MM,F,P)$ is countable and has a single sort descriptor $\SM$; for a countable set of sort descriptors, the same arguments apply with one copy of the set of natural numbers for each sort descriptor.
We write $\Sx_j$ for the variable $\Sx_j^{\SM}$, and, for an atomic formula $\psi$ in $(\Sx_j)_{j<m}$ and a tuple $\bar{k}=(k_j)_{j<m}$ of natural numbers, we write $\psi(\Sx_{\bar{k}})$ for the atomic formula obtained from $\psi$ by substituting $\Sx_{k_j}$ for $\Sx_j$, for every $j<m$; the two formulas have the same predicate symbol, and hence the same a priori bound.

\subsection*{The space of structures on the natural numbers}
Let $\Xs$ be the set of $V$-structures $\fM$ with $\SM^\fM=\NN$.
A structure in $\Xs$ is determined by the interpretations of the symbols of $V$, and $\Xs$ is identified with the product
\[
\prod_{\Sp}[-r_{\Sp},r_{\Sp}]^{\NN^{n_{\Sp}}}\times\prod_{\Sf}\NN^{\NN^{n_{\Sf}}},
\]
over the predicate symbols $\Sp$ in $P$ other than $\unity$ and the function symbols $\Sf$ in $F$, with the product topology and the discrete topology on $\NN$.
As a countable product of Polish spaces, $\Xs$ is a Polish space.
For every formula $\phi(\Sx_0,\dots,\Sx_{m-1})$ of $\Loo^V$ and every tuple $\bar{k}\in\NN^m$, the map $\fM\mapsto\phi^\fM(\bar{k})$ is Borel on $\Xs$: when $\phi$ is atomic, the map is continuous, since, by induction on terms, the value of a term at $\bar{k}$ is locally constant on $\Xs$, the set $\NN$ carrying the discrete topology, so that near each structure the value of $\phi$ at $\bar{k}$ depends on finitely many coordinates; the connectives $\Sneg$, $\Splus$ and $\Sor$ are continuous; and the values of $\supx\phi$ and $\join\Phi$ are suprema of countably many Borel functions.

Let $A$ be the set of pairs $(\psi,\bar{k})$ in which, for some natural number $m$, the formula $\psi$ is an atomic formula in $(\Sx_j)_{j<m}$ with $r_\psi>0$ and $\bar{k}$ belongs to $\NN^m$.
The set $A$ is countable, and the product
\[
Y\coloneq\prod_{(\psi,\bar{k})\in A}[-r_\psi,r_\psi]
\]
is a compact metrizable space.
The \textbf{atomic diagram} of a structure $\fM\in\Xs$ is the point $D(\fM)$ of $Y$ with coordinates $D(\fM)(\psi,\bar{k})\coloneq\psi^\fM(\bar{k})$, and the map $D\colon\Xs\to Y$ is continuous.
For points $y$ and $z$ of $Y$, let
\[
\rho(y,z)\coloneq\sup_{(\psi,\bar{k})\in A}\frac{|y(\psi,\bar{k})-z(\psi,\bar{k})|}{r_\psi},
\]
the uniform distance normalized by the a priori bounds; as a supremum of continuous functions, $\rho$ is lower semicontinuous on $Y\times Y$.

The set $\Surj(\NN)$ of surjections of $\NN$ onto $\NN$ is the intersection, over $n\in\NN$, of the open subsets $\{g\colon n\in g(\NN)\}$ of the Baire space $\NN^\NN$, and is therefore a Polish space~\cite[Theorem~3.11]{Kechris:1995}; under composition, $\Surj(\NN)$ is a monoid.
For $g\in\Surj(\NN)$ and $y\in Y$, let $g^*y$ be the point of $Y$ with coordinates $(g^*y)(\psi,\bar{k})\coloneq y(\psi,g\bar{k})$, where $g\bar{k}=(g(k_j))_{j<m}$.
The map $(g,y)\mapsto g^*y$ is continuous, since each coordinate of $g^*y$ depends on finitely many values of $g$ and on one coordinate of $y$; moreover $(g\circ h)^*=h^*\circ g^*$ for $g,h\in\Surj(\NN)$, and $\rho(g^*y,g^*z)\leq\rho(y,z)$.

A surjection onto the universe of a structure also determines a structure on the natural numbers.
Let $\fM$ be a $V$-structure whose universe $\SM^\fM$ is countable, and let $g\colon\NN\to\SM^\fM$ be a surjection.
The \textbf{pullback} $g^*\fM$ is the member of $\Xs$ in which, for every predicate symbol $\Sp$ and every function symbol $\Sf$ of arity $n$ and every tuple $\bar{k}\in\NN^n$,
\[
\Sp^{g^*\fM}(\bar{k})\coloneq\Sp^\fM(g\bar{k}),\qquad
\Sf^{g^*\fM}(\bar{k})\coloneq\min g^{-1}\bigl(\Sf^\fM(g\bar{k})\bigr),
\]
the least preimage under $g$ of the value of $\Sf^\fM$.
By induction on terms, $g(\tau^{g^*\fM}(\bar{k}))=\tau^\fM(g\bar{k})$ for every term $\tau(\Sx_0,\dots,\Sx_{m-1})$ and every $\bar{k}\in\NN^m$; consequently, $\psi^{g^*\fM}(\bar{k})=\psi^\fM(g\bar{k})$ for every atomic formula $\psi(\Sx_0,\dots,\Sx_{m-1})$, and $D(g^*\fM)=g^*D(\fM)$ when $\fM$ belongs to $\Xs$.
When $V$ has function symbols, the pullbacks chosen by least preimages need not compose as an action of the monoid $\Surj(\NN)$ on $\Xs$; the action is the reindexing of atomic diagrams, and the proofs below use only the identity $D(g^*\fM)=g^*D(\fM)$ and the exact equivalence of $g^*\fM$ with $\fM$, established below.

The comparison of structures uses the distance given by back-and-forth systems.
For $V$-structures $\fM$ and $\fN$, let $\delta(\fM,\fN)$ be the infimum of the real numbers $r>0$ for which there is an $r$-back-and-forth system between $\fM$ and $\fN$ (Definition~\ref{def:system}).
The values of an atomic formula $\psi$ lie in $[-r_\psi,r_\psi]$, so the set of all pairs of tuples of the same length is a $2$-back-and-forth system, and $\delta(\fM,\fN)\leq2$; the structures $\fM$ and $\fN$ are approximately equivalent exactly when $\delta(\fM,\fN)=0$.
An \textbf{exact back-and-forth system} between $\fM$ and $\fN$ is a set of pairs of tuples that contains the pair of empty tuples and satisfies conditions~(ii) and~(iii) of Definition~\ref{def:system} and condition~(i) with $r=0$, that is, with equality of the values of every atomic formula; the structures are \textbf{exactly equivalent} if there is such a system between them.
Exactly equivalent structures are approximately equivalent.
For a surjection $g\colon\NN\to\SM^\fM$, the pairs $(\bar{k},g\bar{k})$, for the tuples $\bar{k}$ of natural numbers, form an exact back-and-forth system between $g^*\fM$ and $\fM$: condition~(ii) holds with the element $g(k)$ for $k\in\NN$, and condition~(iii) with a preimage under $g$ of the given element.
Every $V$-structure whose universe is countable is therefore exactly equivalent to a member of $\Xs$, and the space $\Xs$ represents all such structures.

\begin{lemma}
\label{lem:surjections}
For all structures $\fM$ and $\fN$ in $\Xs$,
\[
\delta(\fM,\fN)=\inf_{p,q\in\Surj(\NN)}\rho\bigl(p^*D(\fM),q^*D(\fN)\bigr),
\]
and $\fM$ and $\fN$ are exactly equivalent if and only if $p^*D(\fM)=q^*D(\fN)$ for some surjections $p$ and $q$.
\end{lemma}

\begin{proof}
Let surjections $p$ and $q$ and a real number $r\geq\rho(p^*D(\fM),q^*D(\fN))$ be given, and let $I$ be the set of pairs $(p\bar{k},q\bar{k})$ for the tuples $\bar{k}$ of natural numbers.
The set $I$ contains the pair of empty tuples.
For condition~(i) of Definition~\ref{def:system}, let $\psi$ be an atomic formula in $(\Sx_j)_{j<m}$ and let $\bar{k}\in\NN^m$; if $r_\psi=0$, both values of $\psi$ vanish, and otherwise $|\psi^\fM(p\bar{k})-\psi^\fN(q\bar{k})|$ is the absolute value of the coordinate $(\psi,\bar{k})$ of $p^*D(\fM)-q^*D(\fN)$, and is at most $r\,r_\psi$.
For condition~(ii), an element $c$ of $\fM$ is $p(j)$ for some $j\in\NN$, since $p$ is surjective, and the pair $(p\bar{k}\,p(j),q\bar{k}\,q(j))$ belongs to $I$; condition~(iii) holds in the same way.
Hence $I$ is an $r$-back-and-forth system for every such $r>0$, and an exact one when $p^*D(\fM)=q^*D(\fN)$; in particular, $\delta(\fM,\fN)$ is at most the infimum in the statement.

Conversely, let $I$ be an $r$-back-and-forth system between $\fM$ and $\fN$, for a real number $r>0$, or an exact one, for $r=0$.
Choosing pairs alternately by conditions~(ii) and~(iii), each extending the previous one, we obtain sequences $(a_j)_{j\in\NN}$ and $(b_j)_{j\in\NN}$ of natural numbers such that the pair $((a_j)_{j<m},(b_j)_{j<m})$ belongs to $I$ for every $m$, and every natural number occurs among the numbers $a_j$ and among the numbers $b_j$.
The maps $p\colon j\mapsto a_j$ and $q\colon j\mapsto b_j$ are surjections.
Let $(\psi,\bar{k})\in A$, with $\psi$ in $(\Sx_j)_{j<m}$, and let $m'$ be a natural number larger than every entry of $\bar{k}$.
The formula $\psi(\Sx_{\bar{k}})$ is an atomic formula in $(\Sx_j)_{j<m'}$ with a priori bound $r_\psi$, and the values of $\psi(\Sx_{\bar{k}})$ at the tuples $(a_j)_{j<m'}$ and $(b_j)_{j<m'}$ are $\psi^\fM(p\bar{k})$ and $\psi^\fN(q\bar{k})$; condition~(i) gives $|\psi^\fM(p\bar{k})-\psi^\fN(q\bar{k})|\leq r\,r_\psi$.
Hence $\rho(p^*D(\fM),q^*D(\fN))\leq r$, so the infimum in the statement is at most $\delta(\fM,\fN)$, and $p^*D(\fM)=q^*D(\fN)$ when $I$ is exact.
\end{proof}

Two structures in $\Xs$ are therefore approximately equivalent when, for every real number $r>0$, the orbits of their atomic diagrams under the monoid $\Surj(\NN)$ come within distance $r$ for $\rho$, and exactly equivalent when the two orbits have a common point.
In classical logic, the group of permutations of $\NN$ plays this part.
Without equality the group does not suffice, since a structure is exactly equivalent to the structure obtained by giving each element a twin.

\begin{example}
\label{ex:twins}
Assume that the only predicate symbol of $V$ besides $\unity$ is $\Sp\colon\SM\times\SM\to[-1,1]$, and that $V$ has no function symbols.
Let $\fM_1$ and $\fM_2$ be the members of $\Xs$ in which $\Sp^{\fM_1}(a,b)$ is $1$ for $a=b$ and $0$ otherwise, and $\Sp^{\fM_2}(a,b)$ is $1$ for $\lfloor a/2\rfloor=\lfloor b/2\rfloor$ and $0$ otherwise.
The surjection $p\colon a\mapsto\lfloor a/2\rfloor$ satisfies $p^*D(\fM_1)=D(\fM_2)$, so that $\fM_1$ and $\fM_2$ are exactly equivalent by Lemma~\ref{lem:surjections}.
For permutations $p$ and $q$ of $\NN$, however, let $a$ and $b$ be the natural numbers with $q(a)=0$ and $q(b)=1$; the coordinate $(\Sp(\Sx_0,\Sx_1),(a,b))$ of $q^*D(\fM_2)$ is $1$, and that of $p^*D(\fM_1)$ is $0$, since $p(a)\neq p(b)$.
Thus $\rho(p^*D(\fM_1),q^*D(\fM_2))=1$ for all permutations $p$ and $q$.
\end{example}

\subsection*{Scott's theorem}
The following theorem is stated in the language of the space $\Xs$ alone.

\begin{theorem}[Scott's theorem in topological form]
\label{thm:topological-scott}
Assume that the vocabulary $V$ is countable and has a single sort descriptor, and let $\Xs$, $Y$, $D$, $\rho$ and the action of $\Surj(\NN)$ on $Y$ be as above.
For $\fM$ and $\fN$ in $\Xs$, let
\[
\delta(\fM,\fN)\coloneq\inf_{p,q\in\Surj(\NN)}\rho\bigl(p^*D(\fM),q^*D(\fN)\bigr).
\]
Then $\delta$ is a pseudometric on $\Xs$, and, for every $\fM\in\Xs$, there is a Borel function $f_\fM\colon\Xs\to[0,2]$, Lipschitz for $\delta$, whose zero set is $\{\fN\in\Xs\colon\delta(\fM,\fN)=0\}$.
In particular, every class of the equivalence relation $\delta(\fM,\fN)=0$ is a Borel subset of $\Xs$, and two members of $\Xs$ lie in the same class if and only if no Borel function on $\Xs$ that is Lipschitz for $\delta$ separates them.
\end{theorem}

\begin{proof}
By Lemma~\ref{lem:surjections}, the function $\delta$ is the distance given by back-and-forth systems.
The pairs $(\bar{k},\bar{k})$ form an exact back-and-forth system between a structure and itself, the inverse of an $r$-back-and-forth system is an $r$-back-and-forth system, and, for an $r$-back-and-forth system $I$ between $\fM$ and $\fN$ and an $s$-back-and-forth system $J$ between $\fN$ and $\fP$, the pairs $(\bar{a},\bar{c})$ for which some tuple $\bar{b}$ satisfies $(\bar{a},\bar{b})\in I$ and $(\bar{b},\bar{c})\in J$ form an $(r+s)$-back-and-forth system between $\fM$ and $\fP$; hence $\delta$ is a pseudometric.

Let $\sigma_\fM$ be the Scott sentence of $\fM$ (Definition~\ref{def:scott-sentence}), and let $f_\fM(\fN)\coloneq\sigma_\fM^\fN$.
The values of $\sigma_\fM$ are nonnegative and the a priori bound of $\sigma_\fM$ is at most $2$, so $f_\fM$ takes values in $[0,2]$, and $f_\fM$ is Borel, as the value of a sentence.
For $\fN$ and $\fN'$ in $\Xs$ and an $r$-back-and-forth system between them, Lemma~\ref{lem:transfer}, applied to the pair of empty tuples, gives $|\sigma_\fM^\fN-\sigma_\fM^{\fN'}|\leq2r$; hence $|f_\fM(\fN)-f_\fM(\fN')|\leq2\,\delta(\fN,\fN')$.
By Theorem~\ref{thm:scott}, $f_\fM(\fN)=0$ exactly when $\fM$ and $\fN$ are approximately equivalent, that is, when $\delta(\fM,\fN)=0$.
The class of $\fM$ is the zero set of the Borel function $f_\fM$, and is therefore Borel.
A function that is Lipschitz for $\delta$ is constant on every class, and $f_\fM$ separates the class of $\fM$ from every other class.
\end{proof}

In classical logic, the theorem gives back the fact that the isomorphism classes of countable structures are Borel.

\begin{corollary}
\label{cor:classical-orbits}
Assume that $V$ is countable and has a single sort descriptor $\SM$ and a predicate symbol $\symb{E}\colon\SM\times\SM\to[-1,1]$.
Let $X_{\symb{E}}$ be the set of members of $\Xs$ in which every atomic formula $\psi$ other than $\unity$ takes only the two values $r_\psi$ and $-r_\psi$ and $\symb{E}$ takes the value $1$ exactly on the diagonal.
For $\fM$ and $\fN$ in $X_{\symb{E}}$, the equality $\delta(\fM,\fN)=0$ holds if and only if $\fN=g^*\fM$ for some permutation $g$ of $\NN$, and every orbit of the group of permutations of $\NN$ on $X_{\symb{E}}$ is a Borel subset of $\Xs$.
\end{corollary}

\begin{proof}
For a permutation $g$, the pullback $g^*\fM$ is the structure transported from $\fM$ along $g$, and $g$ is an isomorphism from $g^*\fM$ onto $\fM$; conversely, an isomorphism $g$ from $\fN$ onto $\fM$ is a permutation with $\fN=g^*\fM$.
By Theorem~\ref{thm:scott} and Corollary~\ref{cor:scott-classical}, $\delta(\fM,\fN)=0$ exactly when $\fN$ is isomorphic to $\fM$.
The set $X_{\symb{E}}$ is closed in $\Xs$, and the orbit of $\fM$ is the intersection of $X_{\symb{E}}$ with the class of $\fM$, a Borel set by Theorem~\ref{thm:topological-scott}.
\end{proof}

\subsection*{López-Escobar's theorem and the limits of the topological form}
Theorem~\ref{thm:topological-scott} separates the classes of approximate equivalence by Borel functions that are Lipschitz for $\delta$.
Theorem~\ref{thm:lopez-escobar} below represents all bounded Borel functions invariant under exact equivalence by the sentences of a logic with continuous connectives, and Example~\ref{ex:exact} shows that some of these functions separate approximately equivalent structures, so that the two classes of functions differ.

Let $\Loo^{c,V}$ be the logic obtained from $\Loo^V$ by adding, for every integer $\ell\geq1$ and every continuous function $u\colon\RR^\ell\to\RR$, a connective that forms from formulas $\phi_1,\dots,\phi_\ell$ the formula $u(\phi_1,\dots,\phi_\ell)$, with value $u(\phi_1^\fM(\bar{a}),\dots,\phi_\ell^\fM(\bar{a}))$ and with a priori bound the maximum of $|u|$ on the product of the intervals $[-r_{\phi_i},r_{\phi_i}]$; joins are formed, as in $\Loo^V$, from sequences of formulas with a common a priori bound and finitely many free variables.
By induction on formulas, the value of every formula of $\Loo^{c,V}$ lies within the a priori bound of the formula, the map $\fM\mapsto\phi^\fM(\bar{k})$ is Borel on $\Xs$ for every formula $\phi$ and every tuple $\bar{k}$, and exactly equivalent structures give every sentence the same value; for the last assertion, the proof of Lemma~\ref{lem:transfer} applies with equalities in place of the inequalities.

\begin{example}
\label{ex:exact}
Assume that the only predicate symbol of $V$ besides $\unity$ is $\Sp\colon\SM\to[-1,1]$, and that $V$ has no function symbols.
Let $\fM$ be a member of $\Xs$ in which $\Sp$ takes each of the values $\frac12$ and $\frac12+\frac1k$, for integers $k\geq2$, at infinitely many elements and takes no other value, and let $\fN$ be a member of $\Xs$ in which $\Sp$ takes each of the values $\frac12+\frac1k$, for integers $k\geq2$, at infinitely many elements and takes no other value.
For a real number $r>0$, the pairs $(\bar{a},\bar{b})$ of tuples of the same length with $|\Sp^\fM(a_i)-\Sp^\fN(b_i)|\leq r$ for every $i$ form an $r$-back-and-forth system, since every value of $\Sp$ in either structure lies within $r$ of the values of $\Sp$ at infinitely many elements of the other; thus $\delta(\fM,\fN)=0$.
For integers $k\geq1$, let $u_k(t)\coloneq\max\{0,1-k|t-\frac12|\}$, and let $\theta$ be the sentence $\synt{\sup_{\Sx_0}}\meet_k u_k(\Sp(\Sx_0))$ of $\Loo^{c,V}$, with a priori bound $1$.
The value of $\theta$ is $1$ in $\fM$ and $0$ in $\fN$.
The function $\fN'\mapsto\theta^{\fN'}$ is therefore Borel on $\Xs$ and invariant under exact equivalence, but not Lipschitz for $\delta$; by Lemma~\ref{lem:transfer}, no sentence of $\Loo^V$ has the same values.
As in Example~\ref{ex:dense}, the obstruction is the absence of a common a priori bound: $\theta$ distinguishes the value $\frac12$ from values arbitrarily close to $\frac12$.
\end{example}

López-Escobar's theorem for $\Loo^{c,V}$ states that the sentences of $\Loo^{c,V}$ give all the invariant Borel functions.
Ben Yaacov, Doucha, Nies and Tsankov proved a version of López-Escobar's theorem for separable metric structures~\cite{BenYaacov-Doucha-Nies-Tsankov:2017}; the version below requires no metric.

\begin{theorem}[López-Escobar's theorem for exact equivalence]
\label{thm:lopez-escobar}
Assume that the vocabulary $V$ is countable and has a single sort descriptor, and let $f\colon\Xs\to\RR$ be a bounded function.
Then $f$ is Borel and takes the same value at exactly equivalent structures if and only if there is a sentence $\phi$ of $\Loo^{c,V}$ with $f(\fM)=\phi^\fM$ for every $\fM\in\Xs$.
\end{theorem}

The proof rests on Vaught's transforms~\cite{Vaught:1974}, \cite[Section~16.B]{Kechris:1995}, taken here over the monoid $\Surj(\NN)$.
For a finite sequence $s=(s_j)_{j<n}$ of natural numbers, let $U_s$ be the set of surjections $g$ with $g(j)=s_j$ for every $j<n$.
The sets $U_s$ are nonempty and open, they form a basis of the topology of $\Surj(\NN)$, and every nonempty open subset of $U_s$ includes a set $U_t$ for some finite sequence $t$ that extends $s$.
Let $F\colon Y\to\RR$ be a bounded Borel function, let $y\in Y$, and let $U$ be a nonempty open subset of $\Surj(\NN)$.
The map $g\mapsto F(g^*y)$ is Borel on $U$, and therefore has the Baire property~\cite[Proposition~8.22]{Kechris:1995}.
The \textbf{Vaught transforms} of $F$ are
\[
\begin{gathered}F^{*U}(y)\coloneq\inf\bigl\{t\in\RR\colon\{g\in U\colon F(g^*y)>t\}\text{ is meager}\bigr\}
\\
F_{*U}(y)\coloneq-(-F)^{*U}(y),\end{gathered}
\]
the essential supremum and the essential infimum, in the sense of category, of the function $g\mapsto F(g^*y)$ on $U$.

\begin{lemma}
\label{lem:localization}
Let $F\colon Y\to\RR$ be a bounded Borel function, let $y\in Y$, and let $s$ be a finite sequence of natural numbers.
Then $F_{*U_s}(y)$ is the infimum of the numbers $F^{*U_t}(y)$ over the finite sequences $t$ that extend $s$.
\end{lemma}

\begin{proof}
Let $h(g)\coloneq F(g^*y)$ and $e\coloneq F_{*U_s}(y)$.
For a finite sequence $t$ that extends $s$ and a real number $\epsilon>0$, the set $\{g\in U_s\colon h(g)>e-\epsilon\}$ is comeager in $U_s$, hence in the open subset $U_t$, and is therefore not meager in $U_t$; hence $F^{*U_t}(y)\geq e-\epsilon$, and so $F^{*U_t}(y)\geq e$.
Conversely, let $q>e$ be a real number.
The set $\{g\in U_s\colon h(g)<q\}$ is not meager and has the Baire property, and is therefore comeager in some nonempty open subset of $U_s$~\cite[Proposition~8.26]{Kechris:1995}, and hence in a set $U_t$ for a finite sequence $t$ that extends $s$.
Then $F^{*U_t}(y)\leq q$.
\end{proof}

A function $F\colon Y\to\RR$ is \textbf{basic} if, for some integer $\ell\geq1$, there are points $(\psi_i,\bar{k}_i)$ of $A$, for $i=1,\dots,\ell$, and a continuous function $u\colon\RR^\ell\to\RR$ such that
\[
F(y)=u\bigl(y(\psi_1,\bar{k}_1),\dots,y(\psi_\ell,\bar{k}_\ell)\bigr)
\]
for every $y\in Y$.

\begin{lemma}
\label{lem:transforms}
Let $c\geq0$ be a real number and let $F\colon Y\to[-c,c]$ be a Borel function.
For every natural number $n$, there are formulas $\phi^+_{F,n}$ and $\phi^-_{F,n}$ of $\Loo^{c,V}$, with free variables among $\Sx_0,\dots,\Sx_{n-1}$ and with a priori bounds at most $c+2$, such that
\[
F^{*U_s}(D(\fM))=\bigl(\phi^+_{F,n}\bigr)^\fM(s)
\qquad\text{and}\qquad
F_{*U_s}(D(\fM))=\bigl(\phi^-_{F,n}\bigr)^\fM(s)
\]
for every $\fM\in\Xs$ and every $s\in\NN^n$.
\end{lemma}

\begin{proof}
Let $\cC$ be the class of Borel functions $G\colon Y\to[-c-2,c+2]$ for which such formulas exist for every $n$, with a priori bounds at most $c+2$.

Let $G$ be a basic function with values in $[-c-2,c+2]$, given by points $(\psi_i,\bar{k}_i)$ of $A$, for $i=1,\dots,\ell$, and a continuous function $u$, and let $J$ be the finite set of the entries of the tuples $\bar{k}_i$.
For $\fM\in\Xs$ and $g\in\Surj(\NN)$, the value $G(g^*D(\fM))=u\bigl(\psi_1^\fM(g\bar{k}_1),\dots,\psi_\ell^\fM(g\bar{k}_\ell)\bigr)$ depends only on the restriction of $g$ to $J$.
For $s\in\NN^n$, the surjections in $U_s$ take the prescribed values $s_j$ at the points $j<n$ of $J$, and the surjections in $U_s$ that extend a given assignment of natural numbers to the points $j\geq n$ of $J$ form a nonempty open set, on which the function $g\mapsto G(g^*D(\fM))$ is constant.
The essential supremum on $U_s$ is therefore the supremum over all such assignments, repeated values included, and is the value at $s$ of the formula
\[
\phi^+_{G,n}\coloneq\synt{\sup_{\Sx_{j_1}}}\cdots\synt{\sup_{\Sx_{j_e}}}\,u\bigl(\psi_1(\Sx_{\bar{k}_1}),\dots,\psi_\ell(\Sx_{\bar{k}_\ell})\bigr),
\]
where $j_1,\dots,j_e$ are the points of $J$ that are at least $n$; the formula $\phi^-_{G,n}$ is the same with infima in place of suprema.
The points $(\psi_i,\bar{k}_i)$ need not be distinct, so the maximum of $|u|$ on the product of the intervals $[-r_{\psi_i},r_{\psi_i}]$ may exceed that of $|G|$; composing $u$ with the retraction of $\RR$ onto $[-c-2,c+2]$, which changes no value of $G$, we may assume that $|u|\leq c+2$ everywhere.
The formula under the quantifiers is the connective $u$ applied to atomic formulas, and the a priori bound of $\phi^\pm_{G,n}$, the maximum of $|u|$ on the product of the intervals $[-r_{\psi_i},r_{\psi_i}]$, is at most $c+2$; so $G\in\cC$.
No condition of distinctness enters the formula; at this step surjections, unlike permutations, dispense with equality.

If $G\in\cC$, then $-G\in\cC$, with $\phi^\pm_{-G,n}\coloneq\synt{-}\phi^\mp_{G,n}$.
Let $G=\sup_kG_k$ for a sequence $(G_k)_{k\in\NN}$ in $\cC$.
A countable union of meager sets is meager, so $G^{*U}=\sup_kG_k^{*U}$ for every nonempty open set $U$, and we let $\phi^+_{G,n}\coloneq\join_k\phi^+_{G_k,n}$.
By Lemma~\ref{lem:localization}, $G_{*U_s}(D(\fM))$ is the infimum, over the integers $m\geq n$ and the tuples $t\in\NN^m$ that extend $s$, of the values $G^{*U_t}(D(\fM))$, and we let
\[
\phi^-_{G,n}\coloneq\meet_{m\geq n}\Sinf{\Sx_n}\cdots\Sinf{\Sx_{m-1}}\phi^+_{G,m}.
\]
Hence $\cC$ is closed under countable suprema and, with negation, under countable infima, and therefore under pointwise limits of sequences, since $\lim_kG_k=\inf_j\sup_{k\geq j}G_k$; every a priori bound in these formulas is at most $c+2$.

The basic functions form an algebra of continuous functions on the compact metrizable space $Y$ that contains the constants and separates points, and the Stone--Weierstrass theorem gives, for every continuous function $G\colon Y\to[-c,c]$, basic functions $G_k$ with $|G-G_k|\leq2^{-k}$ for $k\in\NN$; the functions $G_k+2^{-k}$ are basic, with values in $[-c,c+2]$, and $G$ is their infimum.
Hence $\cC$ contains every continuous function from $Y$ to $[-c,c]$.
Composing with the retraction of $\RR$ onto $[-c,c]$, a continuous function, we see that every Borel function from $Y$ to $[-c,c]$ is obtained from continuous functions from $Y$ to $[-c,c]$ by iterated pointwise limits of sequences~\cite[Theorem~24.3]{Kechris:1995}, and therefore belongs to $\cC$.
\end{proof}

\begin{proof}[Proof of Theorem~\ref{thm:lopez-escobar}]
The values of a sentence of $\Loo^{c,V}$ form a bounded Borel function on $\Xs$ that takes the same value at exactly equivalent structures.

Conversely, let $f$ be such a function, with values in $[-c,c]$.
Two members of $\Xs$ with the same atomic diagram are exactly equivalent, by Lemma~\ref{lem:surjections}, so $f=F\circ D$ for a function $F$ on the set $D(\Xs)$, which is analytic, as the continuous image of a Polish space.
For a rational number $q$, the sets $F^{-1}((q,\infty))=D(f^{-1}((q,\infty)))$ and $D(\Xs)\setminus F^{-1}((q,\infty))=D(f^{-1}((-\infty,q]))$ are disjoint analytic sets, and Lusin's separation theorem~\cite[Theorem~14.7]{Kechris:1995} gives a Borel subset $B_q$ of $Y$ with $F^{-1}((q,\infty))=B_q\cap D(\Xs)$.
Let $B'_q$ be the union of the sets $B_{q'}$ for the rational numbers $q'\geq q$, and, for $y\in Y$, let $G(y)$ be the supremum of the rational numbers $q\in[-c,c]$ with $y\in B'_q$, or $-c$ if there are none.
For $-c\le t<c$, the set $\{y\colon G(y)>t\}$ is the union of the Borel sets $B'_q$ for the rational numbers $q\in(t,c]$;
for $t<-c$ it is $Y$, and for $t\ge c$ it is empty.
Thus $G\colon Y\to[-c,c]$ is Borel.
For $y\in D(\Xs)$, the point $y$ belongs to $B'_q$ exactly when $F(y)>q$, so $G(y)=F(y)$.
Let $\fM\in\Xs$.
For every $g\in\Surj(\NN)$, the structure $g^*\fM$ is exactly equivalent to $\fM$, and
\[
G(g^*D(\fM))=G(D(g^*\fM))=f(g^*\fM)=f(\fM).
\]
The function $g\mapsto G(g^*D(\fM))$ is constant on $\Surj(\NN)=U_\varnothing$, for the empty sequence $\varnothing$, and Lemma~\ref{lem:transforms} gives
\[
f(\fM)=G^{*U_\varnothing}(D(\fM))=\bigl(\phi^+_{G,0}\bigr)^\fM.
\]
The sentence $\phi^+_{G,0}$ of $\Loo^{c,V}$ is the required one.
\end{proof}

Theorems~\ref{thm:topological-scott} and~\ref{thm:lopez-escobar} together mark the reach of the topological form.
The functions of Theorem~\ref{thm:topological-scott} are the values of sentences of $\Loo^V$, which are Lipschitz for $\delta$ by Lemma~\ref{lem:transfer}, and hence constant on the classes of approximate equivalence.
By Theorem~\ref{thm:lopez-escobar}, the Borel functions invariant under exact equivalence are the values of the sentences of $\Loo^{c,V}$, and Example~\ref{ex:exact} shows that some of them separate approximately equivalent structures.
To replace the separating functions by all bounded Borel functions invariant under exact equivalence is therefore to change the equivalence that they detect, from approximate equivalence to exact equivalence.
Whether the sentences of $\Loo^V$ give all the bounded Borel functions that are Lipschitz for $\delta$ is the remaining question; when the atomic values range over finite sets with a uniform gap, they do.

\begin{theorem}[Uniformly separated finite values: representation in the original syntax]
\label{thm:finite-values}
Assume that the vocabulary $V$ is countable and has a single sort descriptor.
For every predicate symbol $\Sp$ of $V$ other than $\unity$, let $\Lambda_{\Sp}$ be a finite subset of $[-r_{\Sp},r_{\Sp}]$, and assume that there is a real number $\gamma>0$ with $|\lambda-\lambda'|>\gamma\,r_{\Sp}$ for every predicate symbol $\Sp$ other than $\unity$ and all distinct $\lambda,\lambda'\in\Lambda_{\Sp}$.
Let $X_\Lambda$ be the set of members of $\Xs$ in which every predicate symbol $\Sp$ other than $\unity$ takes values in $\Lambda_{\Sp}$.
Then a bounded function $f\colon X_\Lambda\to\RR$ is Borel and takes the same value at approximately equivalent structures if and only if there is a sentence $\phi$ of $\Loo^V$ with $f(\fM)=\phi^\fM$ for every $\fM\in X_\Lambda$.
\end{theorem}

\begin{proof}
Sentences of $\Loo^V$ give Borel functions that take the same value at approximately equivalent structures, by Corollary~\ref{cor:equivalence}.
For the converse, we may assume that $X_\Lambda$ is nonempty, since otherwise every sentence represents $f$.
For an atomic formula $\psi$ with predicate symbol $\Sp$ other than $\unity$, let $\Lambda_\psi\coloneq\Lambda_{\Sp}$, and let $\Lambda_{\unity}\coloneq\{1\}$.
For $r<\gamma$, condition~(i) of Definition~\ref{def:system} between members of $X_\Lambda$ forces equality of the values of every atomic formula, so an $r$-back-and-forth system between members of $X_\Lambda$ is exact, and $f$ takes the same value at exactly equivalent members of $X_\Lambda$.
The closed subset $Y_\Lambda\coloneq\prod_{(\psi,\bar{k})\in A}\Lambda_\psi$ of $Y$ is compact metrizable, invariant under the action of $\Surj(\NN)$, and includes $D(X_\Lambda)$; the set $X_\Lambda$ is closed in $\Xs$ and invariant under pullbacks by surjections.
The proofs of Lemmas~\ref{lem:localization} and~\ref{lem:transforms} and of Theorem~\ref{thm:lopez-escobar} apply with $X_\Lambda$ and $Y_\Lambda$ in place of $\Xs$ and $Y$, and yield formulas of $\Loo^V$, once the basic functions are given by formulas of $\Loo^V$ with a priori bounds that depend only on $c$ and $\gamma$.

Let $u\colon\RR^\ell\to\RR$ be a continuous function, and let $\psi_1,\dots,\psi_\ell$ be atomic formulas with $r_{\psi_i}>0$.
On $Y_\Lambda$ the arguments of $u$ range over the finite set $\Lambda_u\coloneq\prod_{i=1}^\ell\Lambda_{\psi_i}$.
Choose a rational number $C$ with $\max_{\Lambda_u}|u|\leq C\leq\max_{\Lambda_u}|u|+1$ and, for each $i$, a rational number $c_i$ with $1/(\gamma r_{\psi_i})\leq c_i\leq2/(\gamma r_{\psi_i})$, and, for $\lambda=(\lambda_i)_{i=1}^\ell\in\Lambda_u$, let
\[
\chi_\lambda\coloneq\unity\Sdot\bigvee_{i=1}^\ell\synt{c_i\cdot}\synt{|}\psi_i\Splus\synt{-}\underline{\lambda_i}\,\synt{|},
\qquad
\eta\coloneq\bigvee_{\lambda\in\Lambda_u}\Bigl(\underline{u(\lambda)}\Sand\synt{C\cdot}\bigl(\synt{2\cdot}\chi_\lambda\Splus\synt{-}\unity\bigr)\Bigr),
\]
where the finite joins are iterations of $\Sor$ and the constants are given by Lemma~\ref{lem:constants} with $\epsilon=r_{\psi_i}$ for $\underline{\lambda_i}$ and $\epsilon=1$ for $\underline{u(\lambda)}$.
At a point of $\Lambda_u$, the formula $\chi_\lambda$ takes the value $1$ if the point is $\lambda$, and the value $0$ otherwise, since some coordinate then differs from $\lambda_i$ by more than $\gamma\,r_{\psi_i}$; so the member of the join indexed by $\lambda$ takes the value $u(\lambda)$ at $\lambda$ and the value $-C\leq u$ elsewhere, and $\eta$ takes the value of $u$ at every point of $\Lambda_u$.
The a priori bound of $\chi_\lambda$ is at most $1+\max_ic_i\cdot3r_{\psi_i}\leq1+6/\gamma$, and that of $\eta$ is at most $C(3+12/\gamma)+1$.
In the proof of Lemma~\ref{lem:transforms}, after the retraction of $u$ onto $[-c-2,c+2]$, the maximum of $|u|$ is at most $c+2$, so that $C\leq c+3$ for the basic functions that occur, and the a priori bounds of all the formulas in the proof are therefore at most a number that depends only on $c$ and $\gamma$.
\end{proof}

For structures in which every atomic formula $\psi$ other than $\unity$ takes only the values $r_\psi$ and $-r_\psi$, Theorem~\ref{thm:finite-values} is López-Escobar's theorem for classical $\Loo$ without equality, with the monoid of surjections in place of the group of permutations.
In general, the question is the following.

\begin{question}
\label{q:lipschitz-sentences}
Assume that the vocabulary $V$ is countable and has a single sort descriptor.
Is every bounded Borel function on $\Xs$ that is Lipschitz for $\delta$ given by a sentence of $\Loo^V$?
\end{question}

A positive answer to Question~\ref{q:lipschitz-sentences} would complement the Scott theorems with a representation theorem in the original syntax, as López-Escobar's theorem complements Scott's theorem in classical logic.
Theorem~\ref{thm:finite-values} is the positive case: when the atomic values range over finite sets with a uniform gap $\gamma$, distinct classes of approximate equivalence lie at distance at least $\gamma$ for $\delta$, so that every bounded function constant on the classes is Lipschitz on $X_\Lambda$; when such a function is also Borel, Theorem~\ref{thm:finite-values} supplies a sentence, with an a priori bound that depends only on the range of the function and on $\gamma$.
When the atomic values accumulate, the problem is to represent these Lipschitz functions while keeping the a priori bounds uniform through the construction.
Example~\ref{ex:exact} shows why invariance under exact equivalence alone does not suffice; its separating function is not Lipschitz for $\delta$, so it leaves the representation problem open.

\end{document}